\documentclass[12pt]{article}
\usepackage{amssymb,amsthm,amsmath,amsfonts,latexsym,tikz,hyperref}
\usepackage[hmargin=1in,vmargin=1in]{geometry}

\newtheorem{thm}{Theorem}[section]
\newtheorem{prop}[thm]{Proposition}
\newtheorem{cor}[thm]{Corollary}
\newtheorem{lem}[thm]{Lemma}
\newtheorem{conj}[thm]{Conjecture}
\newtheorem{exa}[thm]{Example}

\DeclareMathOperator{\Cat}{Cat}

\newcommand{\B}{B}

\newcommand{\ben}{\begin{enumerate}}
\newcommand{\een}{\end{enumerate}}
\newcommand{\ble}{\begin{lem}}
\newcommand{\ele}{\end{lem}}
\newcommand{\bth}{\begin{thm}}
\renewcommand{\eth}{\end{thm}}
\newcommand{\bpr}{\begin{prop}}
\newcommand{\epr}{\end{prop}}
\newcommand{\bco}{\begin{cor}}
\newcommand{\eco}{\end{cor}}
\newcommand{\bcon}{\begin{conj}}
\newcommand{\econ}{\end{conj}}
\newcommand{\bde}{\begin{defn}}
\newcommand{\ede}{\end{defn}}
\newcommand{\bex}{\begin{exa}}
\newcommand{\eex}{\end{exa}}
\newcommand{\barr}{\begin{array}}
\newcommand{\earr}{\end{array}}
\newcommand{\btab}{\begin{tabular}}
\newcommand{\etab}{\end{tabular}}
\newcommand{\beq}{\begin{equation}}
\newcommand{\eeq}{\end{equation}}
\newcommand{\bea}{\begin{eqnarray*}}
\newcommand{\eea}{\end{eqnarray*}}
\newcommand{\bal}{\begin{align*}}
\newcommand{\bce}{\begin{center}}
\newcommand{\ece}{\end{center}}
\newcommand{\bpi}{\begin{picture}}
\newcommand{\epi}{\end{picture}}
\newcommand{\bpp}{\begin{picture}}
\newcommand{\epp}{\end{picture}}
\newcommand{\bfi}{\begin{figure} \begin{center}}
\newcommand{\efi}{\end{center} \end{figure}}
\newcommand{\bprf}{\begin{proof}}
\newcommand{\eprf}{\end{proof}\medskip}
\newcommand{\capt}{\caption}
\newcommand{\bsl}{\begin{slide}{}}
\newcommand{\esl}{\end{slide}}
\newcommand{\bfr}{\begin{frame}}
\newcommand{\efr}{\end{frame}}

\newcommand{\hqed}{\hfill \qed}

\newcommand{\eqqed}[1]{$\rule{1ex}{0ex}\hfill{\dil#1}\hfill\qed$}

\newcommand{\hs}[1]{\hspace{#1}}
\newcommand{\hso}[1]{\hspace{-1pt}}
\newcommand{\vs}[1]{\vspace{#1}}

\newcommand{\emp}{\emptyset}

\newcommand{\sbe}{\subseteq}

\newcommand{\fl}[1]{\lfloor #1 \rfloor}
\newcommand{\ce}[1]{\lceil #1 \rceil}

\def\<{\langle}
\def\>{\rangle}

\newcommand{\ree}[1]{(\ref{#1})}

\newcommand{\be}{\beta}

\newcommand{\de}{\delta}

\newcommand{\io}{\iota}

\newcommand{\la}{\lambda}

\renewcommand{\th}{\theta}

\newcommand{\bbN}{{\mathbb N}}

\newcommand{\cB}{{\cal B}}

\newcommand{\cC}{{\cal C}}

\newcommand{\cS}{{\cal S}}
\newcommand{\cT}{{\cal T}}

\DeclareMathOperator{\wt}{wt}

\newcommand{\dil}{\displaystyle}

\begin{document}
\pagestyle{plain}

\title{Combinatorial interpretations of Lucas analogues of binomial coefficients and Catalan numbers
}
\author{
Curtis Bennett\\[-5pt]
\small Department of Mathematics, California State University,\\[-5pt]
\small Long Beach, CA 90840, USA, {\tt Curtis.Bennett@csulb.edu}\\
Juan Carrillo\\[-5pt]
\small 20707 Berendo Avenue\\[-5pt]
\small Torrance, CA 90502, USA, {\tt juanscarrillo23@gmail.com}\\
John Machacek\\[-5pt]
\small Department of Mathematics and Statistics, York University,\\[-5pt]
\small Toronto, ON M3J 1P3, Canada, {\tt jmachacek.math@gmail.com}\\
Bruce E. Sagan\\[-5pt]
\small Department of Mathematics, Michigan State University,\\[-5pt]
\small East Lansing, MI 48824, USA, {\tt sagan@math.msu.edu}
}

\date{\today\\[10pt]
	\begin{flushleft}
	\end{flushleft}}

\maketitle

\begin{abstract}
The Lucas sequence is a sequence of polynomials in $s,t$ defined recursively by $\{0\}=0$, $\{1\}=1$, and $\{n\}=s\{n-1\}+t\{n-2\}$ for $n\ge2$.  
On specialization of $s$ and $t$ one can recover the Fibonacci numbers, the nonnegative integers, and the $q$-integers $[n]_q$.
Given a quantity which is expressed in terms of products and quotients of nonnegative integers, one obtains a Lucas analogue by replacing each factor of $n$ in the expression with $\{n\}$.  It is then natural to ask if the resulting rational function is actually a polynomial in $s,t$ with nonnegative integer coefficients and, if so, what it counts.  The first simple combinatorial interpretation for this polynomial analogue of the  binomial coefficients was given by Sagan and Savage, although their model resisted being used to prove identities for these Lucasnomials or extending their ideas to other combinatorial sequences.  The purpose of this paper is to give a new, even more natural  model for these Lucasnomials using lattice paths which can be used to prove various equalities as well as extending to Catalan numbers and their relatives, such as those for finite  Coxeter groups.
\end{abstract}

%
%

\section{Introduction}

Let $s$ and $t$ be two indeterminants.  The corresponding {\em Lucas sequence} is defined inductively by letting  $\{0\}=0$, $\{1\}=1$, and
$$
\{n\}=s\{n-1\}+t\{n-2\}
$$
for $n\ge2$.  For example
$$
\{2\}=s, \{3\}=s^2+t, \{4\}=s^3+2st,
$$
and so forth.  Clearly when $s=t=1$ one recovers the Fibonacci sequence. When $s=2$ and $t=-1$ we have $\{n\}=n$.  Furthermore if $s=1+q$ and $t=-q$ then $\{n\}=[n]_q$ where $[n]_q=1+q+\dots+q^{n-1}$ is the usual $q$-analogue of $n$.  So when proving theorems about the Lucas sequence, one gets results about the Fibonacci numbers, the nonnegative integers, 
and $q$-analogues for free.

It is easy to give a combinatorial interpretation to $\{n\}$ in terms of tilings.  Given a row of $n$ squares, let $\cT(n)$ denote the set of tilings $T$ of this strip by monominoes which cover one square and dominoes which cover two adjacent squares.  Figure~\ref{T(3)} shows the tilings in $\cT(3)$.  Given any configuration of tiles $T$ we define its {\em weight} to be
$$
\wt T = s^{\text{number of monominoes in $T$}}\ t^{\text{number of dominoes in $T$}}.
$$
Similarly, given any set of tilings $\cT$ we define its {\em weight} to be 
$$
\wt\cT=\sum_{T\in\cT} \wt T.
$$
To illustrate $\wt(\cT(3))=s^3+2st=\{4\}$.
This presages the next result which follows quickly by induction and is well known so we omit the proof.
\bpr
\label{T(n-1)}
For all $n\ge1$ we have

\vs{3pt}

\eqqed{
\{n\}=\wt(\cT(n-1)).
}
\epr

\bfi
\begin{tikzpicture}
\draw (0,0) grid (3,1);
\fill (.5,.5) circle (.1);
\fill (1.5,.5) circle (.1);
\fill (2.5,.5) circle (.1);
\end{tikzpicture}
\qquad
\begin{tikzpicture}
\draw (0,0) grid (3,1);
\draw (.5,.5)--(1.5,.5);
\fill (.5,.5) circle (.1);
\fill (1.5,.5) circle (.1);
\fill (2.5,.5) circle (.1);
\end{tikzpicture}
\qquad
\begin{tikzpicture}
\draw (0,0) grid (3,1);
\draw (1.5,.5)--(2.5,.5);
\fill (.5,.5) circle (.1);
\fill (1.5,.5) circle (.1);
\fill (2.5,.5) circle (.1);
\end{tikzpicture}
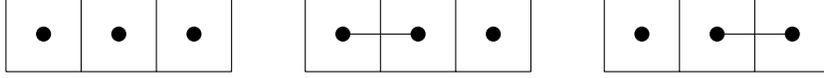
\capt{The tilings in $\cT(3)$}
\label{T(3)}
\efi

Given any quantity which is defined using products and quotients of integers, we can replace each occurrence of $n$ in the expression by $\{n\}$ to obtain its {\em Lucas analogue}.  One can then ask if the resulting rational fiunction is actually a polynomial in $s,t$ with nonnegative integer coefficients and, if it is, whether it is the generating function for some set of combinatorial objects.  
Let $\bbN$ denote the nonnegative integers so that we are interested in showing that various  polynomials are in $\bbN[s,t]$.
We begin by discussing the case of binomial coefficients.

For $n\ge0$, the Lucas analogue of a factorial is the {\em Lucastorial}
$$
\{n\}! = \{1\}\{2\}\dots\{n\}.
$$
Now given $0\le k\le n$ we define the corresponding {\em Lucasnomial} to be
\beq
\label{Lnomial}
{ n \brace k} = \frac{\{n\}!}{\{k\}!\{n-k\}!}.
\eeq
It is not hard to see that this function satisfies an analogue of the binomial recursion (Proposition~\ref{LnomRr} below) and so inductively prove that it is in $\bbN[s,t]$.  

The first simple combinatorial interpretation of the Lucasnomials was given by Sagan and Savage~\cite{ss:cib} using tilings of Young diagrams inside a rectangle.   Earlier but more complicated models were given by Gessel and Viennot~\cite{gv:bdp} and by Benjamin and Plott~\cite{bp:caf}.  Despite its simplicity, there were three difficulties with the Sagan-Savage approach.  The model was not flexible enough to permit combinatorial demonstrations of straight-forward identities involving the Lucasnomials.  Their ideas did not seem to extend to any other related combinatorial sequences such as the Catalan numbers.  And their model contained certain dominoes in the tilings which appeared in an unintuitive manner.  The goal of this paper is to present a new construction which addresses these problems.

We should also mention related work on a $q$-version of these ideas.  As noted above, letting $s=t=1$ reduces $\{n\}$ to $F_n$, the $n$th Fibonacci numbers.  One can then make a $q$-Fibonacci analogue of a quotient of products by replacing each factor of $n$ by $[F_n]_q$.  Working in this framework, some results parallel to ours were found indpendently during the working sessions of the Algebraic
Combinatorics Seminar at the Fields Institute with the active participation of Farid Aliniaeifard, Nantel  Bergeron, Cesar Ceballos, Tom Denton, and Shu Xiao Li~\cite{abcdl}.

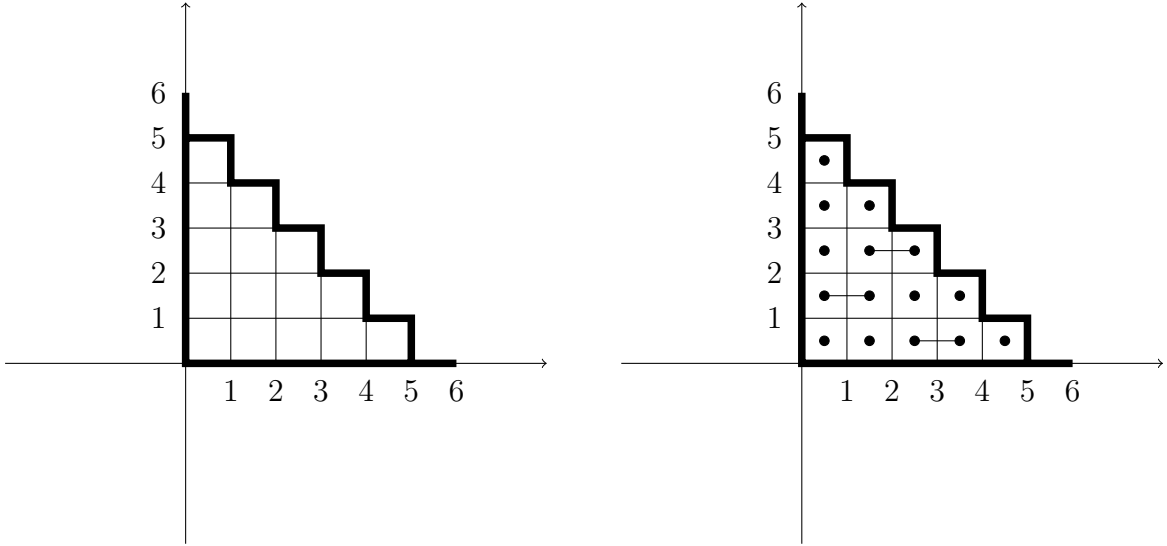
\begin{figure}
\bce
\begin{tikzpicture}[scale=0.6]
\draw (0,6)--(0,0);
\draw (1,5)--(1,0);
\draw (2,4)--(2,0);
\draw (3,3)--(3,0);
\draw (4,2)--(4,0);
\draw (5,1)--(5,0);
\draw (5,0)--(0,0);
\draw (5,1)--(0,1);
\draw (4,2)--(0,2);
\draw (3,3)--(0,3);
\draw (2,4)--(0,4);
\draw (1,5)--(0,5);
\draw[->] (0,-4)--(0,8);
\draw[->] (-4,0)--(8,0);
\draw[line width=1mm] (0,6)--(0,0)--(6,0) (5,0)--(5,1)--(4,1)--(4,2)--(3,2)--(3,3)--(2,3)--(2,4)--(1,4)--(1,5)--(0,5);
\node at (1,-0.6) {$1$};
\node at (2,-0.6) {$2$};
\node at (3,-0.6) {$3$};
\node at (4,-0.6) {$4$};
\node at (5,-0.6) {$5$};
\node at (6,-0.6) {$6$};
\node at (-0.6,1) {$1$};
\node at (-0.6,2) {$2$};
\node at (-0.6,3) {$3$};
\node at (-0.6,4) {$4$};
\node at (-0.6,5) {$5$};
\node at (-0.6,6) {$6$};
\end{tikzpicture}
\hs{20pt}
\begin{tikzpicture}[scale=0.6]
\draw (0,6)--(0,0);
\draw (1,5)--(1,0);
\draw (2,4)--(2,0);
\draw (3,3)--(3,0);
\draw (4,2)--(4,0);
\draw (5,1)--(5,0);
\draw (5,0)--(0,0);
\draw (5,1)--(0,1);
\draw (4,2)--(0,2);
\draw (3,3)--(0,3);
\draw (2,4)--(0,4);
\draw (1,5)--(0,5);
\draw[->] (0,-4)--(0,8);
\draw[->] (-4,0)--(8,0);
\draw[line width=1mm] (0,6)--(0,0)--(6,0) (5,0)--(5,1)--(4,1)--(4,2)--(3,2)--(3,3)--(2,3)--(2,4)--(1,4)--(1,5)--(0,5);
\foreach \x in {.5,1.5,...,4.5 } \filldraw (\x,.5) circle(3pt);
\foreach \x in {.5,1.5,...,3.5 } \filldraw (\x,1.5) circle(3pt);
\foreach \x in {.5,1.5,...,2.5 } \filldraw (\x,2.5) circle(3pt);
\foreach \x in {.5,1.5} \filldraw (\x,3.5) circle(3pt);
\foreach \x in {.5} \filldraw (\x,4.5) circle(3pt);
\draw (2.5,.5)--(3.5,.5) (.5,1.5)--(1.5,1.5) (1.5,2.5)--(2.5,2.5);
\node at (1,-0.6) {$1$};
\node at (2,-0.6) {$2$};
\node at (3,-0.6) {$3$};
\node at (4,-0.6) {$4$};
\node at (5,-0.6) {$5$};
\node at (6,-0.6) {$6$};
\node at (-0.6,1) {$1$};
\node at (-0.6,2) {$2$};
\node at (-0.6,3) {$3$};
\node at (-0.6,4) {$4$};
\node at (-0.6,5) {$5$};
\node at (-0.6,6) {$6$};
\end{tikzpicture}
\ece
\caption{$\de_6$ embedded in $\mathbb{R}^2$ on the left and a tiling on the right.}
\label{R2}
\end{figure}

We will need to consider lattice paths inside tilings of Young diagrams.  Let $\la=(\la_1,\la_2,\dots,\la_l)$ be an integer partition, that is, a weakly decreasing sequence of positive integers.  The $\la_i$ are called {\em parts} and the {\em length} of $\la$ is the number of parts and denoted $l(\la)$.  The {\em Young diagram} of $\la$ is an array of left-justified rows of boxes which we will write in French notation so that $\la_i$ is the number of boxes in the $i$th row from the bottom of the diagram.  We will also use the notation $\la$ for the diagram of $\la$.   Furthermore, we will embed this diagram in the first quadrant of a Cartesian coordinate system with the boxes being unit squares and the southwest-most corner of $\la$ being  the origin.  Finally, it will be convenient in what follows to consider the unit line segments from $(\la_1,0)$ to $(\la_1+1,0)$ and from $(0,l(\la))$ to $(0,l(\la)+1)$ to be part of $\la$'s diagram.  On the left in Figure~\ref{R2} the diagram of $\la=\de_6$ is outlined with thick lines where
$$
\de_n=(n-1,n-2,\dots,1).
$$

A {\em tiling} of $\la$ is a tiling $T$ of the rows of the diagram with monominoes and dominoes.  We let $\cT(\la)$ denote the set of all such tilings.  An element of $\cT(\de_6)$ is shown on the right in Figure~\ref{R2}.  We write $\wt\la$ for the more cumbersome $\wt(\cT(\la))$.  The fact that  $\wt\de_n=\{n\}!$ follows directly from the definitions.  So to prove that $\{n\}!/p(s,t)$ is a polynomial for some polynomial $p(s,t)$, it suffices to partition $\cT(\de_n)$ into subsets, which we will call {\em blocks}, such that $\wt\beta$ is evenly divisible by $p(s,t)$ for all blocks $\beta$.  We will use lattice paths inside $\de_n$ to create the partitions where the choice of path will vary depending on which Lucas analogue we are considering.

The rest of this paper is organized as follows.  In the next section, we give our new combinatorial interpretation for the Lucasnomials.  	In Section~\ref{ifl} we prove two identities using this model.  The demonstration for one of them is straightforward, but the other requires a surprisingly intricate algorithm.  Section~\ref{cfc} is devoted to showing how our model can be modified to give combinatorial interpretations to Lucas analogues of the Catalan and Fuss-Catalan numbers.  In the following section we prove that the Coxeter-Catalan numbers for any Coxeter group have polynomial Lucas analogues and that the same is true for the infiinite families of Coxeter groups in the Fuss-Catalan case.  In fact we generalize these results by considering $d$-divisible diagrams, $d$ being a positive integer, where each row has length one less than a multiple of $d$.  We end with a section containing comments and directions for future research.

\section{Lucasnomials}

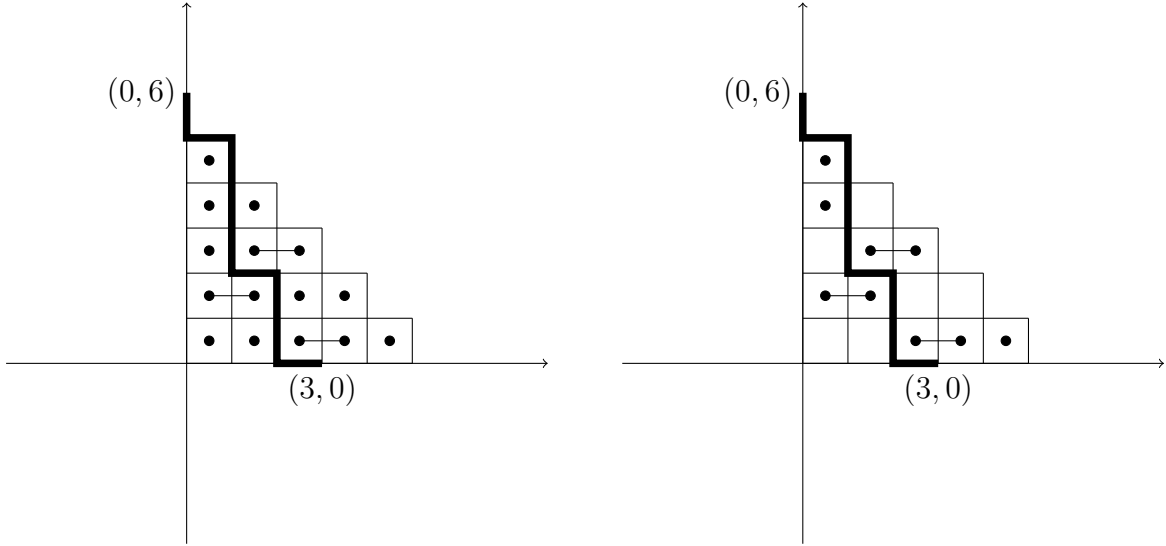
\begin{figure}
\bce
\begin{tikzpicture}[scale=0.6]
\draw (0,6)--(0,0);
\draw (1,5)--(1,0);
\draw (2,4)--(2,0);
\draw (3,3)--(3,0);
\draw (4,2)--(4,0);
\draw (5,1)--(5,0);
\draw (5,0)--(0,0);
\draw (5,1)--(0,1);
\draw (4,2)--(0,2);
\draw (3,3)--(0,3);
\draw (2,4)--(0,4);
\draw (1,5)--(0,5);
\draw[->] (0,-4)--(0,8);
\draw[->] (-4,0)--(8,0);
\draw[line width=1mm] (3,0)--(2,0)--(2,2)--(1,2)--(1,5)--(0,5)--(0,6);
\foreach \x in {.5,1.5,...,4.5 } \filldraw (\x,.5) circle(3pt);
\foreach \x in {.5,1.5,...,3.5 } \filldraw (\x,1.5) circle(3pt);
\foreach \x in {.5,1.5,...,2.5 } \filldraw (\x,2.5) circle(3pt);
\foreach \x in {.5,1.5} \filldraw (\x,3.5) circle(3pt);
\foreach \x in {.5} \filldraw (\x,4.5) circle(3pt);
\draw (2.5,.5)--(3.5,.5) (.5,1.5)--(1.5,1.5) (1.5,2.5)--(2.5,2.5);
\node at (3,-0.6) {$(3,0)$};
\node at (-1,6) {$(0,6)$};
\end{tikzpicture}
\hs{20pt}
\begin{tikzpicture}[scale=0.6]
\draw (0,6)--(0,0);
\draw (1,5)--(1,0);
\draw (2,4)--(2,0);
\draw (3,3)--(3,0);
\draw (4,2)--(4,0);
\draw (5,1)--(5,0);
\draw (5,0)--(0,0);
\draw (5,1)--(0,1);
\draw (4,2)--(0,2);
\draw (3,3)--(0,3);
\draw (2,4)--(0,4);
\draw (1,5)--(0,5);
\draw[->] (0,-4)--(0,8);
\draw[->] (-4,0)--(8,0);
\draw[line width=1mm] (3,0)--(2,0)--(2,2)--(1,2)--(1,5)--(0,5)--(0,6);
\foreach \x in {2.5,3.5,4.5 } \filldraw (\x,.5) circle(3pt);
\foreach \x in {.5,1.5} \filldraw (\x,1.5) circle(3pt);
\foreach \x in {1.5,2.5 } \filldraw (\x,2.5) circle(3pt);
\foreach \x in {.5} \filldraw (\x,3.5) circle(3pt);
\foreach \x in {.5} \filldraw (\x,4.5) circle(3pt);
\draw (2.5,.5)--(3.5,.5) (.5,1.5)--(1.5,1.5) (1.5,2.5)--(2.5,2.5);
\node at (3,-0.6) {$(3,0)$};
\node at (-1,6) {$(0,6)$};
\end{tikzpicture}
\ece
\caption{The path for the  tiling in Figure~\ref{R2} and the corresponding partial tiling}
\label{path}
\end{figure}

In this section we will use the method outlined in the introduction to show that the Lucasnomials defined by~\ree{Lnomial} are polynomials in $s$ and $t$.  In particular, we will prove the following result.
\bth
\label{Lnom:ptn}
Given $0\le k\le n$ there is a partition of $\cT(\de_n)$ such that $\{k\}!\{n-k\}!$ divides $\wt\beta$ for every block $\beta$.
\eth
\bprf
Given $T\in\cT(\de_n)$ we will describe the block $\beta$ containing it by using a lattice path $p$.  The path will start at $(k,0)$ and end at $(0,n)$ taking unit steps north ($N$) and west ($W$).  If $p$ is at a lattice point $(x,y)$ then it moves north to $(x,y+1)$ as long as doing so will not cross a domino and not take $p$ out of the diagram of $\de_n$.  Otherwise, $p$ moves west to $(x-1,y)$.  For example, if $T$ is the tiling in Figure~\ref{R2} then the resulting path is shown on the left in Figure~\ref{path}.  Indeed, initially $p$ is forced west by a domino above, then moves north twice until forced west again by a domino, then moves north three times until doing so again would take it out of $\de_6$, and finishes by following the boundary of the diagram.  In this case we write $p=WNNWNNNWN$.

The north steps of $p$ are of two two kinds: those which are immediately preceded by a west step and those which are not.
Calll the former $NL$ steps (since the $W$ and $N$ step together look like a letter ell) and the latter $NI$ steps.  
In our example the first, third, and sixth north steps are $NL$ while the others are $NI$.
The block containing the original tiling $T$ consists of all tilings in $\cT(\de_n)$ which agree with $T$ to the right of each $NL$ step and to the left of each $NI$ step.  Returning to Figure~\ref{path}, the diagram on the right shows the tiles which are common to all tilings in $T$'s block whereas the squares which are blank can be tiled arbitrarily.    Now there are $k-i$ boxes to the left of the $i$th $NL$ step and thus the weight of tiling these boxes is  $\{k-i+1\}$.  Since there $k$ such steps, the total contribution to $\wt\beta$ is
$\{k\}!$ for the boxes to the left of these steps.  Similarly, $\{n-k\}!$ is the contribution to $\wt\beta$ of the boxes to the right of the $NI$ steps.  This completes the proof.
\eprf

The tiling showing the fixed tiles for a given block $\beta$ in the partition of the previous theorem will be called a {\em binomial partial tiling} $B$.  As just proved, $\wt\beta  = \{k\}!\{n-k\}!\wt B$.  Thus we have the following result.
\bco
\label{bpartial}
Given $0\le k\le n$ we have
$$
{n \brace k}=\sum_B \wt B
$$
where the sum is over all binomial partial tilings associated with lattice paths from $(k,0)$ to $(0,n)$ in $\de_n$.
Thus ${n \brace k}\in\bbN[s,t]$. \hqed
\eco

\begin{figure}
\bce
\begin{tikzpicture}[scale=0.6]
\draw (0,3)--(0,0);
\draw (1,3)--(1,0);
\draw (2,3)--(2,0);
\draw (3,3)--(3,0);
\draw (3,0)--(0,0);
\draw (3,1)--(0,1);
\draw (3,2)--(0,2);
\draw (3,3)--(0,3);
\foreach \x in {.5,1.5,2.5 } \foreach \y in {.5,1.5,2.5 } \filldraw (\x,\y) circle(3pt);
\draw (.5,2.5)--(1.5,2.5) (1.5,.5)--(1.5,1.5) (2.5,.5)--(2.5,1.5);
\draw[line width=1mm] (0,0)--(1,0)--(1,2)--(2,2)--(2,3)--(3,3);
\end{tikzpicture}
\ece
\caption{The tiling of a rectangle  corresponding to the partial tiling in Figure~\ref{path}}
\label{R}
\end{figure}

We end this section by describing the relationship between the tilings we have been considering and those in the model of Sagan and Savage.   In their interpretation, one considered all lattice paths $p$  in a $k\times(n-k)$ rectangle $R$ starting at the southwest corner, ending at the northeast corner, and taking unit steps north and east.  The path divides $R$  into two partitions: $\la$ whose parts are the rows of  boxes in $R$ northwest of $p$ and $\la^*$ whose parts are the columns of $R$ southeast of $p$.  One then considers all tilings of  $R$ which are tilings of $\la$ (so any dominoes are horizontal) and of $\la^*$ (so any dominoes are vertical) such that each tiling of a column of $\la^*$ begins with a domino.  They then proved that ${n\brace k}$ is the generating function for all such tilings.  But there is a bijection between these tilings and our binomial partial tilings where one uses the fixed tilings to left of $NI$ steps for the rows of $\la$ and those to the right of the $NL$ steps for $\la^*$.  Figure~\ref{R} shows the tiling of $R$ corresponding to the partial tiling in Figure~\ref{path}.  Note that dominoes in the tiling of $\la^*$ occur naturally in the context of binomial partial tilings rather than just being an imposed condition.  And, as we will see, the viewpoint of partial tilings is much more flexible than that of tilings of a rectangle.

\section{Identities for Lucasnomials}
\label{ifl}

We will now use Corollary~\ref{bpartial} to prove various identities for Lucasnomials.  We start with the analogue of the binomial recursion mentioned in the introduction.  In the Sagan and Savage paper, this formula was first proved by other means and then used to obtain their combinatorial interpretation.  Here, the recursion follows easily from our model.

\bpr
\label{LnomRr}
For $0<k<n$ we have
$$
{n\brace k} = \{k+1\} {n-1\brace k}  + t \{n-k-1\} {n-1\brace k-1}.
$$
\epr
\bprf
By Corollary~\ref{bpartial} it suffices to partition the set of partial tilings for ${n\brace k}$ into two subsets whose generating functions give the two terms of the recursion.  First consider the binomial partial tilings $B$ whose path $p$ starts from $(k,0)$ with an $N$ step.  Since this is an $NI$ step, the portion of the first row to the left of the step is tiled, and by Proposition~\ref{T(n-1)} the generating function for such tilings is $\{k+1\}$.  Now $p$ continues from $(k,1)$ through the remaining rows which form the partition $\de_{n-1}$.  It follows that the weights of this portion of the corresponding partial tilings sum to ${n-1\brace k}$.  Thus these $p$ give the first term of the recursion. 

Suppose now that $p$ starts with a $W$ step.  It follows that the second step of $p$ must be $N$ and so an $NL$ step.  In this case the portion of the first row to the right of the $NL$ step is tiled and that tiling begins with a domino.  Using the same reasoning as in the previous paragraph, one sees that the weight generating function for the tiling of the first row is $t \{n-k-1\}$ while ${n-1\brace k-1}$ accounts for the rest of the rows of the tiling.  
\eprf

We will next give a bijective proof of the symmetry of the Lucasnomials.  In particular, we will construct an involution to demonstrate the following result.
\begin{prop}
\label{sym}
For  $0\le t\le k\le n$ we have
\beq
\label{symeq}
\{k\}\{k-1\}\dots\{k-t+1\}{n\brace k} = \{n-k+t\}\dots\{n-k+1\}{n \brace n-k+t}.
\eeq
In particular, when $t=0$,
$$
{ n \brace k} = {n\brace n-k}.
$$
\end{prop}

\bfi
\begin{tikzpicture}[scale=0.6]
\draw (0,7)--(0,0);
\draw (1,6)--(1,0);
\draw (2,5)--(2,0);
\draw (3,4)--(3,0);
\draw (4,3)--(4,0);
\draw (5,2)--(5,0);
\draw (6,1)--(6,0);
\draw (7,0)--(7,0);

\draw (7,0)--(0,0);
\draw (6,1)--(0,1);
\draw (5,2)--(0,2);
\draw (4,3)--(0,3);
\draw (3,4)--(0,4);
\draw (2,5)--(0,5);
\draw (1,6)--(0,6);
\draw (0,7)--(0,7);

\draw[->] (0,-1)--(0,8);
\draw[->] (-1,0)--(8,0);

\node at (1,-0.6) {$1$};
\node at (2,-0.6) {$2$};
\node at (3,-0.6) {$3$};
\node at (4,-0.6) {$4$};
\node at (5,-0.6) {$5$};
\node at (6,-0.6) {$6$};
\node at (7,-0.6) {$7$};

\node at (-0.6,1) {$1$};
\node at (-0.6,2) {$2$};
\node at (-0.6,3) {$3$};
\node at (-0.6,4) {$4$};
\node at (-0.6,5) {$5$};
\node at (-0.6,6) {$6$};
\node at (-0.6,7) {$7$};

\foreach \x in {4.5,5.5 } \filldraw (\x,.5) circle(3pt);
\foreach \x in {3.5,4.5 } \filldraw (\x,1.5) circle(3pt);
\foreach \x in {.5,1.5,2.5} \filldraw (\x,2.5) circle(3pt);
\foreach \x in {.5,1.5,2.5} \filldraw (\x,3.5) circle(3pt);
\foreach \x in {5.5,6.5,...,8.5} \filldraw (\x,6.5) circle(3pt);
\foreach \x in {5.5,6.5,7.5} \filldraw (\x,5.5) circle(3pt);

\draw (4.5,.5)--(5.5,.5) 
          (3.5,1.5)--(4.5,1.5) 
          (.5,2.5)--(1.5,2.5)
          (1.5,3.5)--(2.5,3.5);

\draw (6.5,6.5)--(7.5,6.5)
	(5.5,5.5)--(6.5,5.5);

\draw[line width=1mm] (5,0)--(4,0)--(4,1)--(3,1)--(3,2)--(3,3)--(3,4)--(2,4)--(2,5)--(1,5)--(1,6)--(0,6)--(0,7);

\draw (5,5) grid (8,7);
\draw (8,6) grid (9,7);

\node at (-1.6,3.5){$B$};
\node at (4.4,6.5){$S_1$};
\node at (4.4,5.5){$S_2$};

\end{tikzpicture}

\capt{An extended binomial partial tiling of type $(7,5,2)$ \label{752}}
\efi

Although this proposition is easy to prove algebraically, the algorithm giving the involution is surprisingly intricate. To define the bijection we will need the following concepts.
A {\em strip of length $k$} will be a row of $k$ squares, that is, the Young diagram of $\la=(k)$.  For  $0\le t\le k\le n$ an {\em extended binomial partial tiling of type $(n,k,t)$} is a $(t+1)$-tuple $\cB=(B;S_1,\dots,S_t)$ where
\ben
\item $B$ is a partial binomial tiling of $\de_n$ whose lattice path starts at $(k,0)$, and
\item $S_i$ is a tiled strip of length $k-i$ for $1\le i\le t$.
\een
For brevity we will sometimes write $\cB=(B;\cS)$ where $\cS=(S_1,\dots,S_t)$.  In figures, we will display the strips to the northeast of $B$.   See Figure~\ref{752} for an example with $(n,k,t)=(7,5,2)$.  Clearly the sum of the weights of all $\cB$ of type $(n,k,t)$ is the left-hand side of equation~\ree{symeq}.
 Our involution will be a  map $\io:\cB\mapsto\cC$ where $\cB$ and $\cC$  are of types $(n,k,t)$ and $(n,n-k+t,t)$, respectively.  This will provide a combinatorial proof of Proposition~\ref{sym}.

To describe the algorithm producing $\io$ we need certain operations on partial binomial tilings and strips.  Given two strips $R,S$ we denote their concatenation by $RS$ or $R\cdot S$.  So, using the strips in Figure~\ref{752},
\bce
\begin{tikzpicture}[scale=0.6]
\draw (0,0) grid (7,1);
\foreach \x in {.5,1.5,...,6.5} \filldraw (\x,.5) circle(3pt);
\draw (1.5,.5)--(2.5,.5) (4.5,.5)--(5.5,.5);
\node at (-3,.5) {$S_1 S_2=S_1\cdot S_2=$};
\node at (7.5,.2) {.};
\end{tikzpicture}
\ece
Given  a partially tiled strip $R$ of length $n$ and a partial binomial tiling $B$ of $\de_n$  we define their {\em concatenation}, $RB=R\cdot B$, to be the tiling of $\de_{n+1}$ whose first (bottom) row is $R$ and with the remaining rows tiled as in $B$.  So if $B'$ is the partial tiling of $\de_6$ given by the top six rows of the partial tiling $B$ in Figure~\ref{752} then $B=RB'$ where
\bce
\begin{tikzpicture}[scale=0.6]
\draw (0,0) grid (6,1);
\foreach \x in {4.5,5.5} \filldraw (\x,.5) circle(3pt);
\draw (4.5,.5)--(5.5,.5);
\node at (-2,.5) {$R=$};
\node at (6.5,.2) {.};
\end{tikzpicture}
\ece
Note that only for certain $R$ will the concatenation $RB$ remain a partial binomial tiling for some path.  
In particular $R$ will have to be either a {\em left strip} where only the left-most boxes are tiled, or a {\em right strip} with tiles only on the right-most boxes.  The example $R$ is a right strip tiled by a domino.

Given a strip $S$ of length $k$ and $0\le s\le k$ we define $S\fl{s}$ and $S\ce{s}$ to be the strips consisting of, respectively, the first $s$  and the last $s$ boxes of $S$.  Continuing our example
\bce
\begin{tikzpicture}[scale=0.6]
\draw (0,0) grid (3,1);
\foreach \x in {.5,1.5,2.5} \filldraw (\x,.5) circle(3pt);
\draw (1.5,.5)--(2.5,.5);
\node at (-2,.5) {$S_1\fl{3}=$};
\node at (5.5,.5){and $S_1\ce{3}=$};
\draw (8,0) grid (11,1);
\foreach \x in {8.5,9.5,10.5} \filldraw (\x,.5) circle(3pt);
\draw (8.5,.5)--(9.5,.5);
\node at (11.5,.2) {.};
\end{tikzpicture}
\ece
Note that these notations are undefined if taking the desired boxes would involve breaking a domino.  Also, to simplify notation, we will use $RS\fl{s}$ to be the first $s$ boxes of the concatentation of $RS$, while $R\cdot S\fl{s}$ will be the concatenation of $R$ with the first $s$ boxes of $S$.  A similar convention applies to the last $s$ boxes.  We will also use the notation $S^r$ for the {\em reverse} of a strip obtained by reflecting it in a vertical axis.  In our example,
\bce
\begin{tikzpicture}[scale=0.6]
\draw (0,0) grid (3,1);
\foreach \x in {.5,1.5,2.5} \filldraw (\x,.5) circle(3pt);
\draw (1.5,.5)--(2.5,.5);
\node at (-1,.5) {$S_2^r=$};
\node at (3.5,.2) {.};
\end{tikzpicture}
\ece

Our algorithm will break into four cases depending on the following concept.  
Call a point $(r,0)$ an {\em $NI$ point} of a partial binomial tiling $B$  if taking an $N$ step from this vertex stays in $B$ and does not cross a domino.  Otherwise call $(r,0)$ an {\em $NL$ point} which also includes the case where this vertex is not in $B$ to begin with.
We consider use the previous two definitions for strips by considering them as being embedded in the first quadrant as a one-row partition.
Because our algorithm is recursive, we will have to be careful about its notation. 
A priori, given a partial binomial tiling $B$ the notation $\io(B)$ is not well defined since $\io$ needs as input a pair $\cB=(B;\cS)$.  However, it will be convenient to write
$$
\io(B;S_1,\dots,S_t) = (\io(B);\io(S_1),\dots,\io(S_t))
$$
where it is understood on the right-hand side that $\io$ is always taken with respect to the input pair $(B;\cS)$ to the algorithm.  
Because the algorithm is recursive, we will also have to apply $\io$ to $\cB'=(B';S_1',\dots,S_r')$ where $B'$ is $B$ with its bottom row removed and $S_1',\dots,S_r'$ are certain strips.  So we define $\io'(B')$ and $\io'(S_i')$  for $1\le i\le r$ by
$$
\io(B';S'_1,\dots,S'_r) = (\io'(B');\io'(S'_1),\dots,\io'(S'_r)).
$$
In particular, if $S_i'=S_j$ for some $i,j$ then $\io'(S_j)=\io'(S'_i)$ so that $S_j$ is being treated as an element of $\cB'$ rather than of $\cB$.  We now have all the necessary concepts to present the recursive algorithm which is given in Figure~\ref{alg}.  

\bfi
\ben
\item[]{\bf Algorithm $\io$}
\item[] Input: An extended binomial tiling $\cB=(B;S_1,\dots,S_t)$ having type $(n,k,t)$.
\item[] Output: An extended binomial tiling $\io(\cB)=(C;T_1,\dots,T_t)=(C;\cT)$ having type $(n,n-k+t,t)$.
\item If $n=0$ then $\io$ is the identity and $\cC=\cB$.
\item If $n>0$ then let $R$ be the strip of tiled squares in the bottom row of $B$, and let $B'$ be $B$ with the bottom row removed.
\item Construct $\cB'=(B',\cS')$, calculate $\io'(\cB')$ recusively, and then define $\io(\cB)$ using the following four cases. 
\ben
\item  If $(k,0)$ is an $NI$ point of $B$ and $(k-t-1,0)$ is an $NI$ point of $B$ then let
$$
\barr{rcl}
S_{t+1} &=& R\fl{k-t-1},\\ 
\cS' &=&(S_1,\dots,S_t,S_{t+1}),\\
C&= & R^L\cdot\io'(B')  \text{ where $R^L$ is a left strip tiled by $R'=\io'(S_{t+1}) \cdot R\ce{t+1}^r$},\\
\cT&=& (\io'(S_1),\dots,\io'(S_t)).
\earr
$$
\item If $(k,0)$ is an $NI$ point of $B$ and $(k-t-1,0)$ is an $NL$ point of $B$ then let
$$
\barr{rcl}
\cS' &=&(S_1,\dots,S_t),\\
C&= & R^R\cdot\io'(B')  \text{ where $R^R$ is a right strip tiled by $R'=R\fl{k-t}^r$},\\
\cT&=& (\io'(S_t)\cdot R\ce{t}^r,\ \io'(S_1),\ \dots,\ \io'(S_{t-1})).
\earr
$$
\item If $(k,0)$ is an $NL$ point of $B$ and $(k-t-1,0)$ is an $NI$ point of $S_1$ (by convention, this is considered to be true if $t=0$ so that $S_1$ does not exist) then let
$$
\barr{rcl}
S_{t+1} &=& S_1\fl{k-t-1},\\
\cS' &=&(S_2,\dots,S_t,S_{t+1}),\\
C&= & R^L\cdot\io'(B')  \text{ where $R^L$ is a left strip tiled by $R'=R^r S_1^r\fl{n-k+t}$},\\
\cT&=& (\io'(S_2), \dots,\io'(S_t), \io'(S_{t+1})).
\earr
$$
\item  If $(k,0)$ is an $NL$ point of $B$ and $(k-t-1,0)$ is an $NL$ point of $S_1$ then let
$$
\barr{rcl}
\cS' &=&(S_2,\dots,S_t),\\
C&= & R^R\cdot\io'(B')  \text{ where $R^R$ is a right strip tiled by $R'=R^r S_1^r\ce{k-t}$},\\
\cT&=& (R^r S_1^r\fl{n-k+t-1},\ \io'(S_2),\ \dots,\ \io'(S_t)).
\earr
$$
\een
\een
\capt{The algorithm for computing the involution $\io$} \label{alg}
\efi

A step-by-step example of applying $\io$ to the extended binomial tiling in Figure~\ref{752} is begun in Figure~\ref{752I} and finished in Figure~\ref{752II}.  In it, $(B^{(i)};\cS^{(i)})$ represents the pair on which $\io$ is called in the $i$th iteration.  So in the notation of the algorithm $(B^{(1)};\cS^{(1)})=(B';\cS')$ and so forth.  These superscripts will make it clear that when we write, for example, $\io(B^{(i)})$ we are referring to $\io$ acting on the pair $(B^{(i)};\cS^{(i)})$. These pairs are listed down the left sides of the figures.  On the right sides are as much of the output $(C,\cT)$ as has been constructed after each iteration.  Strips contributing to $\cT$ may be left partially blank if the recursion has not gone deep enough yet to completely fill them.  
The circles used for the tiles have been replaced by numbers or letters to make it easier to follow the movement of the tiles.
Tiles from $B$ are labeled with the number of their row while tiles from the strips are labeled alphabetically.
Finally, the ``maps to" symbols indicate which of the four cases (a)--(d) of the algorithm is being used at each step.  We will now consider the first four steps in detail since they will illustrate each of the four cases.  The reader should find it easy to fill in the particulars for the rest of the steps.

\bfi
\vspace*{-30pt}
\begin{tikzpicture}[scale=0.6]
\draw (0,7)--(0,0);
\draw (1,6)--(1,0);
\draw (2,5)--(2,0);
\draw (3,4)--(3,0);
\draw (4,3)--(4,0);
\draw (5,2)--(5,0);
\draw (6,1)--(6,0);
\draw (7,0)--(7,0);

\draw (7,0)--(0,0);
\draw (6,1)--(0,1);
\draw (5,2)--(0,2);
\draw (4,3)--(0,3);
\draw (3,4)--(0,4);
\draw (2,5)--(0,5);
\draw (1,6)--(0,6);
\draw (0,7)--(0,7);

\draw[->] (0,-1)--(0,8);
\draw[->] (-1,0)--(8,0);

\node at (1,-0.6) {$1$};
\node at (2,-0.6) {$2$};
\node at (3,-0.6) {$3$};
\node at (4,-0.6) {$4$};
\node at (5,-0.6) {$5$};
\node at (6,-0.6) {$6$};
\node at (7,-0.6) {$7$};

\node at (-0.6,1) {$1$};
\node at (-0.6,2) {$2$};
\node at (-0.6,3) {$3$};
\node at (-0.6,4) {$4$};
\node at (-0.6,5) {$5$};
\node at (-0.6,6) {$6$};
\node at (-0.6,7) {$7$};

\foreach \x in {4.5,5.5 } \filldraw (\x,.5) node{$1$};
\foreach \x in {3.5,4.5 } \filldraw (\x,1.5) node{$2$};
\foreach \x in {.5,1.5,2.5} \filldraw (\x,2.5) node{$3$};
\foreach \x in {.5,1.5,2.5} \filldraw (\x,3.5) node{$4$};
\foreach \x in {5.5,6.5,...,8.5} \filldraw (\x,6.5) node{$a$};
\foreach \x in {5.5,6.5,7.5} \filldraw (\x,5.5) node{$b$};

\draw (4.6,.5)--(5.4,.5) 
          (3.6,1.5)--(4.4,1.5) 
          (.6,2.5)--(1.4,2.5)
          (1.6,3.5)--(2.4,3.5);

\draw (6.6,6.5)--(7.4,6.5)
	(5.6,5.5)--(6.4,5.5);

\draw[line width=1mm] (5,0)--(4,0)--(4,1)--(3,1)--(3,2)--(3,3)--(3,4)--(2,4)--(2,5)--(1,5)--(1,6)--(0,6)--(0,7);

\draw (5,5) grid (8,7);
\draw (8,6) grid (9,7);

\node at (-1.6,3.5){$B$};
\node at (4.4,6.5){$S_1$};
\node at (4.4,5.5){$S_2$};

\node at (15,3.5){$(C,\cT)=\emp$};

\node at (20,3.5){$\barr{c}\rm (d)\\  \mapsto\earr$};

\end{tikzpicture}
\vspace*{10pt}

\begin{tikzpicture}[scale=0.6]

\draw (0,7)--(0,1);
\draw (1,6)--(1,1);
\draw (2,5)--(2,1);
\draw (3,4)--(3,1);
\draw (4,3)--(4,1);
\draw (5,2)--(5,1);
\draw (6,1)--(6,1);

\draw (6,1)--(0,1);
\draw (5,2)--(0,2);
\draw (4,3)--(0,3);
\draw (3,4)--(0,4);
\draw (2,5)--(0,5);
\draw (1,6)--(0,6);
\draw (0,7)--(0,7);

\draw[->] (0,0)--(0,8);
\draw[->] (-1,1)--(7,1);

\node at (1,.6) {$1$};
\node at (2,.6) {$2$};
\node at (3,.6) {$3$};
\node at (4,.6) {$4$};
\node at (5,.6) {$5$};
\node at (6,.6) {$6$};

\node at (-0.6,2) {$1$};
\node at (-0.6,3) {$2$};
\node at (-0.6,4) {$3$};
\node at (-0.6,5) {$4$};
\node at (-0.6,6) {$5$};
\node at (-0.6,7) {$6$};

\foreach \x in {3.5,4.5 } \filldraw (\x,1.5) node{$2$};
\foreach \x in {.5,1.5,2.5} \filldraw (\x,2.5) node{$3$};
\foreach \x in {.5,1.5,2.5} \filldraw (\x,3.5) node{$4$};

\foreach \x in {5.5,6.5,7.5} \filldraw (\x,6.5) node{$b$};
\draw  (5,6) grid (8,7);

\draw (3.6,1.5)--(4.4,1.5) 
          (.6,2.5)--(1.4,2.5)
          (1.6,3.5)--(2.4,3.5);
          
\draw (5.6,6.5)--(6.4,6.5);

\node at (-1.6,4){$B^{(1)}$};
\node at (4.2,6.5){$S^{(1)}_1$};

\draw[line width=1mm] (4,1)--(3,1)--(3,2)--(3,3)--(3,4)--(2,4)--(2,5)--(1,5)--(1,6)--(0,6)--(0,7);

\node at (8.6,3.5) {$C$};

\draw[->] (9.4,0)--(18,0);
\draw[->] (10,-.6)--(10,7.5);
\draw (10,1)--(16,1)
          (11,1)--(11,0) (12,1)--(12,0) (13,1)--(13,0) (14,1)--(14,0) (15,1)--(15,0) (16,1)--(16,0);
          
\node at (9.4,1){$1$};

\node at (11,-0.6) {$1$};
\node at (12,-0.6) {$2$};
\node at (13,-0.6) {$3$};
\node at (14,-0.6) {$4$};
\node at (15,-0.6) {$5$};
\node at (16,-0.6) {$6$};
\node at (17,-0.6) {$7$};

\foreach \x in {13.5,...,15.5 } \filldraw (\x,.5) node{$a$};

\draw  (13.6,.5)--(14.4,.5);

\draw (15,5) grid (17,7);
\draw (17,6) grid (18,7);

\foreach \x in {15.5,16.5} \filldraw (\x,6.5) node{$1$};
\foreach \x in {17.5 } \filldraw (\x,6.5) node{$a$};

\draw[line width=1mm] (14,0)--(13,0)--(13,1);

\draw (15.6,6.5)--(16.4,6.5);

\node at (14.2,6.5){$T_1$};
\node at (13.7,5.5){$\io(S_1^{(1)})$};

\node at (20,3.5){$\barr{c}\rm (c)\\  \mapsto\earr$};
\end{tikzpicture}

\vspace*{10pt}

\begin{tikzpicture}[scale=0.6]

\node at (-1.6,4.5){$\B^{(2)}$};

\draw (0,7)--(0,2);
\draw (1,6)--(1,2);
\draw (2,5)--(2,2);
\draw (3,4)--(3,2);
\draw (4,3)--(4,2);
\draw (5,2)--(5,2);

\draw (5,2)--(0,2);
\draw (4,3)--(0,3);
\draw (3,4)--(0,4);
\draw (2,5)--(0,5);
\draw (1,6)--(0,6);
\draw (0,7)--(0,7);

\draw[->] (0,1)--(0,8);
\draw[->] (-1,2)--(6,2);

\node at (1,1.6) {$1$};
\node at (2,1.6) {$2$};
\node at (3,1.6) {$3$};
\node at (4,1.6) {$4$};
\node at (5,1.6) {$5$};

\node at (-0.6,3) {$1$};
\node at (-0.6,4) {$2$};
\node at (-0.6,5) {$3$};
\node at (-0.6,6) {$4$};
\node at (-0.6,7) {$5$};

\foreach \x in {.5,1.5,2.5} \filldraw (\x,2.5) node{$3$};
\foreach \x in {.5,1.5,2.5} \filldraw (\x,3.5) node{$4$};

\foreach \x in {5.5,6.5} \filldraw (\x,6.5) node{$b$};
\draw  (5,6) grid (7,7);
\draw (5.6,6.5)--(6.4,6.5);
\node at (4.2,6.5){$S_1^{(2)}$};

\draw  (.6,2.5)--(1.4,2.5)
          (1.6,3.5)--(2.4,3.5);

\draw[line width=1mm] (3,2)--(3,3)--(3,4)--(2,4)--(2,5)--(1,5)--(1,6)--(0,6)--(0,7);

\node at (8.6,3.5) {$C$};

\draw[->] (9.4,0)--(18,0);
\draw[->] (10,-.6)--(10,7.5);
\draw (10,1)--(16,1) (10,2)--(15,2)
          (11,2)--(11,0) (12,2)--(12,0) (13,2)--(13,0) (14,2)--(14,0) (15,2)--(15,0) (16,1)--(16,0);
          
\node at (9.4,1){$1$};
\node at (9.4,2){$2$};

\node at (11,-0.6) {$1$};
\node at (12,-0.6) {$2$};
\node at (13,-0.6) {$3$};
\node at (14,-0.6) {$4$};
\node at (15,-0.6) {$5$};
\node at (16,-0.6) {$6$};
\node at (17,-0.6) {$7$};

\foreach \x in {13.5,...,15.5 } \filldraw (\x,.5) node{$a$};
\foreach \x in {10.5,11.5 } \filldraw (\x,1.5) node{$2$};
\foreach \x in {12.5} \filldraw (\x,1.5) node{$b$};

\draw  (13.6,.5)--(14.4,.5) (10.6,1.5)--(11.4,1.5);

\draw (15,5) grid (17,7);
\draw (17,6) grid (18,7);

\foreach \x in {15.5,16.5} \filldraw (\x,6.5) node{$1$};
\foreach \x in {17.5 } \filldraw (\x,6.5) node{$a$};

\draw (15.6,6.5)--(16.4,6.5);

\node at (14.2,6.5){$T_1$};
\node at (13.7,5.5){$\io(S_1^{(2)})$};

\draw[line width=1mm] (14,0)--(13,0)--(13,1)--(13,2);

\node at (20,3.5){$\barr{c}\rm (b)\\  \mapsto\earr$};
\end{tikzpicture}

\vspace*{10pt}

\begin{tikzpicture}[scale=0.6]

\node at (-1.6,5){$\B^{(3)}$};

\draw (0,7)--(0,3);
\draw (1,6)--(1,3);
\draw (2,5)--(2,3);
\draw (3,4)--(3,3);
\draw (4,3)--(4,3);

\draw (4,3)--(0,3);
\draw (3,4)--(0,4);
\draw (2,5)--(0,5);
\draw (1,6)--(0,6);
\draw (0,7)--(0,7);

\draw[->] (0,2)--(0,8);
\draw[->] (-1,3)--(5,3);

\node at (1,2.6) {$1$};
\node at (2,2.6) {$2$};
\node at (3,2.6) {$3$};
\node at (4,2.6) {$4$};

\node at (-0.6,4) {$1$};
\node at (-0.6,5) {$2$};
\node at (-0.6,6) {$3$};
\node at (-0.6,7) {$4$};

\foreach \x in {.5,1.5,2.5} \filldraw (\x,3.5) node{$4$};

\foreach \x in {5.5,6.5} \filldraw (\x,6.5) node{$b$};
\draw  (5,6) grid (7,7);
\draw (5.6,6.5)--(6.4,6.5);
\node at (4.2,6.5){$S_1^{(3)}$};

\draw (1.6,3.5)--(2.4,3.5);

\draw[line width=1mm] (3,3)--(3,4)--(2,4)--(2,5)--(1,5)--(1,6)--(0,6)--(0,7);

\node at (8.6,4){$C$};
\draw[->] (9.4,0)--(18,0);
\draw[->] (10,-.6)--(10,7.5);
\draw (10,1)--(16,1) (10,2)--(15,2) (10,3)--(14,3)
          (11,3)--(11,0) (12,3)--(12,0) (13,3)--(13,0) (14,3)--(14,0) (15,2)--(15,0) (16,1)--(16,0);
          
\node at (9.4,1){$1$};
\node at (9.4,2){$2$};
\node at (9.4,3){$3$};

\node at (11,-0.6) {$1$};
\node at (12,-0.6) {$2$};
\node at (13,-0.6) {$3$};
\node at (14,-0.6) {$4$};
\node at (15,-0.6) {$5$};
\node at (16,-0.6) {$6$};
\node at (17,-0.6) {$7$};

\foreach \x in {13.5,...,15.5 } \filldraw (\x,.5) node{$a$};
\foreach \x in {10.5,11.5 } \filldraw (\x,1.5) node{$2$};
\foreach \x in {12.5} \filldraw (\x,1.5) node{$b$};
\foreach \x in {12.5,13.5} \filldraw (\x,2.5) node{$3$};

\draw (13.6,.5)--(14.4,.5) (10.6,1.5)--(11.4,1.5) (12.6,2.5)--(13.4,2.5);

\draw (15,5) grid (17,7);
\draw (17,6) grid (18,7);

\foreach \x in {15.5,16.5} \filldraw (\x,6.5) node{$1$};
\foreach \x in {17.5 } \filldraw (\x,6.5) node{$a$};
\filldraw (16.5,5.5) node{$3$};

\draw (15.6,6.5)--(16.4,6.5);

\node at (14.2,6.5){$T_1$};
\node at (13.7,5.5){$\io(S_1^{(2)})$};
\draw[->] (15.5,4.3)--(15.5,4.9);
\node at (15.5,3.9){$\io(S_1^{(3)})$};

\draw[line width=1mm] (14,0)--(13,0)--(13,1)--(13,2)--(13,2)--(12,2)--(12,3);

\node at (20,3.5){$\barr{c}\rm (a)\\  \mapsto\earr$};

\end{tikzpicture}

\capt{The $\io$ recursion applied to the extended binomial partial tiling in Figure~\ref{752}, part 1} \label{752I}
\efi

\bfi
\vspace*{-30pt}
\begin{tikzpicture}[scale=0.6]

\node at (-1.6,5.5){$\B^{(4)}$};

\draw (0,7)--(0,4);
\draw (1,6)--(1,4);
\draw (2,5)--(2,4);
\draw (3,4)--(3,4);

\draw (3,4)--(0,4);
\draw (2,5)--(0,5);
\draw (1,6)--(0,6);
\draw (0,7)--(0,7);

\draw[->] (0,3)--(0,8);
\draw[->] (-1,4)--(5,4);

\node at (1,3.6) {$1$};
\node at (2,3.6) {$2$};
\node at (3,3.6) {$3$};

\node at (-0.6,5) {$1$};
\node at (-0.6,6) {$2$};
\node at (-0.6,7) {$3$};

\foreach \x in {5.5,6.5} \filldraw (\x,6.5) node{$b$};
\filldraw (5.5,5.5) node{$4$};
\draw  (5,6) grid (7,7);
\draw (5,5) grid (6,6);
\draw (5.6,6.5)--(6.4,6.5);
\node at (4.2,6.5){$S_1^{(4)}$};
\node at (4.2,5.5){$S_2^{(4)}$};

\draw[line width=1mm] (3,4)--(2,4)--(2,5)--(1,5)--(1,6)--(0,6)--(0,7);

\node at (8.6,4){$C$};
\draw[->] (9.4,0)--(18,0);
\draw[->] (10,-.6)--(10,7.5);
\draw (10,1)--(16,1) (10,2)--(15,2) (10,3)--(14,3) (10,4)--(13,4)
          (11,4)--(11,0) (12,4)--(12,0) (13,4)--(13,0) (14,3)--(14,0) (15,2)--(15,0) (16,1)--(16,0);
          
\node at (9.4,1){$1$};
\node at (9.4,2){$2$};
\node at (9.4,3){$3$};
\node at (9.4,4){$4$};

\node at (11,-0.6) {$1$};
\node at (12,-0.6) {$2$};
\node at (13,-0.6) {$3$};
\node at (14,-0.6) {$4$};
\node at (15,-0.6) {$5$};
\node at (16,-0.6) {$6$};
\node at (17,-0.6) {$7$};

\foreach \x in {13.5,...,15.5 } \filldraw (\x,.5) node{$a$};
\foreach \x in {10.5,11.5 } \filldraw (\x,1.5) node{$2$};
\foreach \x in {12.5} \filldraw (\x,1.5) node{$b$};
\foreach \x in {12.5,13.5} \filldraw (\x,2.5) node{$3$};
\foreach \x in {10.5,11.5}\filldraw (\x,3.5) node{$4$};

\draw (13.6,.5)--(14.4,.5) (10.6,1.5)--(11.4,1.5) (12.6,2.5)--(13.4,2.5) (10.6,3.5)--(11.4,3.5);

\draw (15,5) grid (17,7);
\draw (17,6) grid (18,7);

\foreach \x in {15.5,16.5} \filldraw (\x,6.5) node{$1$};
\foreach \x in {17.5 } \filldraw (\x,6.5) node{$a$};
\filldraw (16.5,5.5) node{$3$};

\draw (15.6,6.5)--(16.4,6.5);

\node at (14.2,6.5){$T_1$};
\node at (13.7,5.5){$\io(S_1^{(2)})$};
\draw[->] (15.5,4.3)--(15.5,4.9);
\node at (15.5,3.9){$\io(S_1^{(4)})$};

\draw[line width=1mm] (14,0)--(13,0)--(13,1)--(13,2)--(13,2)--(12,2)--(12,3)--(12,4);

\node at (20,3.5){$\barr{c}\rm (c)\\  \mapsto\earr$};

\end{tikzpicture}

\vspace*{10pt}

\begin{tikzpicture}[scale=0.6]

\node at (-1.6,6){$\B^{(5)}$};

\draw (0,7)--(0,5);
\draw (1,6)--(1,5);
\draw (2,5)--(2,5);

\draw (2,5)--(0,5);
\draw (1,6)--(0,6);
\draw (0,7)--(0,7);

\draw[->] (0,5)--(0,8);
\draw[->] (-1,5)--(3,5);

\node at (1,4.6) {$1$};
\node at (2,4.6) {$2$};

\node at (-0.6,6) {$1$};
\node at (-0.6,7) {$2$};

\filldraw (5.5,6.5) node{$4$};
\draw  (5,6) grid (6,7);
\draw (5,5)--(5,6);
\node at (4.2,6.5){$S_1^{(5)}$};
\node at (4.2,5.5){$S_2^{(5)}$};

\draw[line width=1mm] (2,5)--(1,5)--(1,6)--(0,6)--(0,7);

\node at (8.6,4){$C$};
\draw[->] (9.4,0)--(18,0);
\draw[->] (10,-.6)--(10,7.5);
\draw (10,1)--(16,1) (10,2)--(15,2) (10,3)--(14,3) (10,4)--(13,4) (10,5)--(12,5)
          (11,5)--(11,0) (12,5)--(12,0) (13,4)--(13,0) (14,3)--(14,0) (15,2)--(15,0) (16,1)--(16,0);
          
\node at (9.4,1){$1$};
\node at (9.4,2){$2$};
\node at (9.4,3){$3$};
\node at (9.4,4){$4$};
\node at (9.4,5){$5$};

\node at (11,-0.6) {$1$};
\node at (12,-0.6) {$2$};
\node at (13,-0.6) {$3$};
\node at (14,-0.6) {$4$};
\node at (15,-0.6) {$5$};
\node at (16,-0.6) {$6$};
\node at (17,-0.6) {$7$};

\foreach \x in {13.5,...,15.5 } \filldraw (\x,.5) node{$a$};
\foreach \x in {10.5,11.5 } \filldraw (\x,1.5) node{$2$};
\foreach \x in {12.5} \filldraw (\x,1.5) node{$b$};
\foreach \x in {12.5,13.5} \filldraw (\x,2.5) node{$3$};
\foreach \x in {10.5,11.5}\filldraw (\x,3.5) node{$4$};
\foreach \x in {10.5,11.5}\filldraw (\x,4.5) node{$b$};

\draw (13.6,.5)--(14.4,.5) (10.6,1.5)--(11.4,1.5) (12.6,2.5)--(13.4,2.5) (10.6,3.5)--(11.4,3.5)  (10.6,4.5)--(11.4,4.5);

\draw (15,5) grid (17,7);
\draw (17,6) grid (18,7);

\foreach \x in {15.5,16.5} \filldraw (\x,6.5) node{$1$};
\foreach \x in {17.5 } \filldraw (\x,6.5) node{$a$};
\filldraw (16.5,5.5) node{$3$};

\draw (15.6,6.5)--(16.4,6.5);

\node at (14.2,6.5){$T_1$};
\node at (13.7,5.5){$\io(S_1^{(2)})$};
\draw[->] (15.5,4.3)--(15.5,4.9);
\node at (15.5,3.9){$\io(S_1^{(5)})$};

\draw[line width=1mm] (14,0)--(13,0)--(13,1)--(13,2)--(13,2)--(12,2)--(12,3)--(12,4)--(12,5);
\node at (20,3.5){$\barr{c}\rm (d)\\  \mapsto\earr$};
\end{tikzpicture}
\vspace*{10pt}

\begin{tikzpicture}[scale=0.6]

\node at (-1.6,6.6){$\B^{(6)}$};

\draw (0,7)--(0,6);
\draw (1,6)--(1,6);

\draw (1,6)--(0,6);
\draw (0,7)--(0,7);

\draw[->] (0,5)--(0,8);
\draw[->] (-1,6)--(2,6);

\node at (1,5.6) {$1$};

\node at (-0.6,7) {$1$};

\draw[line width=1mm] (1,6)--(0,6)--(0,7);

\draw (5,6)--(5,7);
\node at (4.2,6.5){$S_1^{(6)}$};

\node at (8.6,4){$C$};
\draw[->] (9.4,0)--(18,0);
\draw[->] (10,-.6)--(10,7.5);
\draw (10,1)--(16,1) (10,2)--(15,2) (10,3)--(14,3) (10,4)--(13,4) (10,5)--(12,5) (10,6)--(11,6)
          (11,6)--(11,0) (12,5)--(12,0) (13,4)--(13,0) (14,3)--(14,0) (15,2)--(15,0) (16,1)--(16,0);
          
\node at (9.4,1){$1$};
\node at (9.4,2){$2$};
\node at (9.4,3){$3$};
\node at (9.4,4){$4$};
\node at (9.4,5){$5$};
\node at (9.4,6){$6$};
\node at (9.4,7){$7$};

\node at (11,-0.6) {$1$};
\node at (12,-0.6) {$2$};
\node at (13,-0.6) {$3$};
\node at (14,-0.6) {$4$};
\node at (15,-0.6) {$5$};
\node at (16,-0.6) {$6$};
\node at (17,-0.6) {$7$};

\foreach \x in {13.5,...,15.5 } \filldraw (\x,.5) node{$a$};
\foreach \x in {10.5,11.5 } \filldraw (\x,1.5) node{$2$};
\foreach \x in {12.5} \filldraw (\x,1.5) node{$b$};
\foreach \x in {12.5,13.5} \filldraw (\x,2.5) node{$3$};
\foreach \x in {10.5,11.5}\filldraw (\x,3.5) node{$4$};
\foreach \x in {10.5,11.5}\filldraw (\x,4.5) node{$b$};

\draw (13.6,.5)--(14.4,.5) (10.6,1.5)--(11.4,1.5) (12.6,2.5)--(13.4,2.5) (10.6,3.5)--(11.4,3.5)  (10.6,4.5)--(11.4,4.5);

\draw (15,5) grid (17,7);
\draw (17,6) grid (18,7);

\foreach \x in {15.5,16.5} \filldraw (\x,6.5) node{$1$};
\foreach \x in {17.5 } \filldraw (\x,6.5) node{$a$};
\draw (16.5,5.5) node{$3$};
\draw (15.5,5.5) node{$4$};

\draw (15.6,6.5)--(16.4,6.5);

\node at (14.2,6.5){$T_1$};
\node at (13.7,5.5){$\io(S_1^{(2)})$};

\draw[line width=1mm] (14,0)--(13,0)--(13,1)--(13,2)--(13,2)--(12,2)--(12,3)--(12,4)--(12,5)--(11,5)--(11,6);
\node at (20,3.5){$\barr{c}\rm (d)\\  \mapsto\earr$};
\end{tikzpicture}

\vspace*{10pt}

\begin{tikzpicture}[scale=0.6]

\node at (-1.6,7){$\B^{(7)}$};

\draw[->] (0,6)--(0,8);
\draw[->] (-1,7)--(1,7);
\filldraw (0,7) circle(3pt);

\node at (8.6,4){$C$};
\draw[->] (9.4,0)--(18,0);
\draw[->] (10,-.6)--(10,7.5);
\draw (10,1)--(16,1) (10,2)--(15,2) (10,3)--(14,3) (10,4)--(13,4) (10,5)--(12,5) (10,6)--(11,6)
          (11,6)--(11,0) (12,5)--(12,0) (13,4)--(13,0) (14,3)--(14,0) (15,2)--(15,0) (16,1)--(16,0);
          
\node at (9.4,1){$1$};
\node at (9.4,2){$2$};
\node at (9.4,3){$3$};
\node at (9.4,4){$4$};
\node at (9.4,5){$5$};
\node at (9.4,6){$6$};
\node at (9.4,7){$7$};

\node at (11,-0.6) {$1$};
\node at (12,-0.6) {$2$};
\node at (13,-0.6) {$3$};
\node at (14,-0.6) {$4$};
\node at (15,-0.6) {$5$};
\node at (16,-0.6) {$6$};
\node at (17,-0.6) {$7$};

\foreach \x in {13.5,...,15.5 } \filldraw (\x,.5) node{$a$};
\foreach \x in {10.5,11.5 } \filldraw (\x,1.5) node{$2$};
\foreach \x in {12.5} \filldraw (\x,1.5) node{$b$};
\foreach \x in {12.5,13.5} \filldraw (\x,2.5) node{$3$};
\foreach \x in {10.5,11.5}\filldraw (\x,3.5) node{$4$};
\foreach \x in {10.5,11.5}\filldraw (\x,4.5) node{$b$};

\draw (13.6,.5)--(14.4,.5) (10.6,1.5)--(11.4,1.5) (12.6,2.5)--(13.4,2.5) (10.6,3.5)--(11.4,3.5)  (10.6,4.5)--(11.4,4.5);

\draw (15,5) grid (17,7);
\draw (17,6) grid (18,7);

\foreach \x in {15.5,16.5} \filldraw (\x,6.5) node{$1$};
\foreach \x in {17.5 } \filldraw (\x,6.5) node{$a$};
\draw (16.5,5.5) node{$3$};
\draw (15.5,5.5) node{$4$};

\draw (15.6,6.5)--(16.4,6.5);

\node at (14.2,6.5){$T_1$};
\node at (14.2,5.5){$T_2$};

\draw[line width=1mm] (14,0)--(13,0)--(13,1)--(13,2)--(13,2)--(12,2)--(12,3)--(12,4)--(12,5)--(11,5)--(11,6)--(10,6)--(10,7);
\end{tikzpicture}
\hspace*{50pt}
\capt{The $\io$ recursion applied to the extended binomial partial tiling in Figure~\ref{752}, part 2} \label{752II}
\efi

Initially $(n,k,t)=(7,5,2)$.  We see that $(5,0)$ is an $NL$ point of $B$ and $(5-2-1,0)=(2,0)$ is also an $NL$ point of $S_1$.  So we are in case (d).  Accordingly, $\cS^{(1)}=(S_2)$ so that $S_1^{(1)}=S_2$ which is the strip filled with $b$'s.  Also
\bce
\begin{tikzpicture}[scale=0.6]
\draw (0,0) grid (6,1);
\foreach\x in {.5,1.5} \draw (\x,.5) node{$1$};
\foreach \x in {2.5,3.5,4.5,5.5} \draw (\x,.5) node{$a$};
\draw (.6,.5)--(1.4,.5) (3.6,.5)--(4.4,.5);
\node at (-2,.5) {$R^r S_1^r=$};
\node at (6.5,.2) {.};
\end{tikzpicture}
\ece
Taking that last $5-2=3$ squares of $R^r S_1^r$ gives a right strip for $C$.  Recalling that $S_2=S_1^{(1)}$ we have that 
$\cT=(T_1,\io(S_1^{(1)}))$ where $T_1$ is the tiling of the remaining $7-5+2-1=3$ squares of $R^r S_1^r$.  Since we will have to recurse further to compute $\io(S_1^{(1)})$, the squares for that strip are left blank in the figure.

In $\cB^{(1)}$ we have a $(6,4,1)$ extended tiling.  Furthermore $(4,0)$ is an $NL$ point of $B^{(1)}$ while $(4-1-1,0)=(2,0)$ is an $NI$ point of $S_1^{(1)}$.  It follows that we are in case (c).  Thus $\cS^{(2)}$ consists of only one strip which are the first two tiles of $S_1^{(1)}$.  The second row of $C$ will be the left strip gotten by taking the tiles in the first  $6-4-1=3$ boxes of the concatenation
\bce
\begin{tikzpicture}[scale=0.6]
\draw (0,0) grid (5,1);
\foreach\x in {.5,1.5} \draw (\x,.5) node{$2$};
\foreach \x in {2.5,3.5,4.5} \draw (\x,.5) node{$b$};
\draw (.6,.5)--(1.4,.5) (3.6,.5)--(4.4,.5);
\node at (5.5,.2) {.};
\end{tikzpicture}
\ece
And $\cT^{(2)}$ consists of the single strip $\io(S_1^{(2)})$ which still remains to be computed.

At the next stage, the extended tiling is of type $(5,3,1)$.  The two points $(3,0)$  and  $(3-1-1,0)=(1,0)$ are, respectively, $NI$ and $NL$ points of $B^{(2)}$.  This is case (b) so $\cS^{(3)}=\cS^{(2)}$.  The new row of $C$ consists of the right strip tiled by the first $5-3=2$ tiles of  the lowest row of $B^{(2)}$ in reverse order.   (Reversal  does nothing since the tiling is just a single domino.)   And the single strip in $\cT^{(3)}$ is obtained by concatenating $\io(S_1^{(3)})$  with the last tile of the lowest row of $B^{(2)}$ in reverse order (which again does nothing since the tiling is just a single monomino).

We have that $\cB^{(3)}$ is of type $(4,3,1)$.  The point $(3,0)$ is an $NI$ step of $B^{(3)}$ as is $(3-1-1,0)=(1,0)$.  So we are in case (a).  So $\cS^{(4)}$ will have a new strip consisting of the first tile of 
\bce
\begin{tikzpicture}[scale=0.6]
\draw (0,0) grid (3,1);
\foreach\x in {.5,1.5,2.5} \draw (\x,.5) node{$4$};
\draw (1.6,.5)--(2.4,.5);
\node at (3.5,.2) {.};
\end{tikzpicture}
\ece
Now $C$ adds a row consisting of the reversal of the tiles on the remaining two squares of the above strip, while $\cT$ does not change from the previous step.

\bth
\label{iothm}
The map $\io$ is a well defined involution on extended binomial tilings.
\eth
\bprf
We induct on $n$ where the case $n=0$ is trivial.  So assume $n>0$ and that the theorem holds for extended binomial tilings with first parameter $n-1$.  We will now go through each of the cases of the algorithm in turn.

Consider case (a).  To check that $\io$ is well defined, we must first show that restricting to the first $k-t-1$ (or the last $t+1$) boxes of $R$ does not break a domino.  But this is true since $R$ has length $k$ and $(k-t-1,0)$ is an $NI$ point of $B$.  
Note also  that $|S_{t+1}|=k-t-1$ is the correct length to be the final strip in $\cS'$ since one takes an $NI$ step to go from $B$ to $B'$ and so the path still has $x$-coordinate $k$.  Similarly, the other strips of $\cS'$ have the appropriate lengths.
Next we must be sure that the left strip used for the bottom of $C$ will permit the beginning of a path starting at $(n-k+t,0)$ with an $NI$ step.  
First note that  the forms of $B'$ and $\cS'$ show that $\cB'$ has parameters $(n-1,k,t+1)$,  so by induction $\io'(\cB')$ has type
$(n-1,n-k+t,t+1)$
It follows that  $\io'(S_{t+1})$ has length $(n-k+t)-(t+1)=n-k-1$, and  the number of boxes tiled in the first row of $C$ is
$$
|\io'(S_{t+1}) \cdot R\ce{t+1}^r| = (n-k-1)+(t+1)=n-k+t
$$
as desired.  We must also make sure that once the $NI$ step is taken in $C$, its end point will be the same as the initial point of the path when we compute $\io'(\cB')$.  But since the first step in $C$ will be $NI$, its $x$-coordinate will still be $n-k+t$ which agrees with the middle parameter computed for $\io'(\cB')$ above.  Finally, we must check that the entries of $\cT$ have the correct lengths.  But this follows from the fact that the middle parameters for $\cB$ and $\cB'$ are both $k$.

We now check that $\io^2(\cB)=\cB$ in case (a).  To avoid confusion we will always use two different alphabets to distinguish between $\cB$ and $\cC=\io(\cB)$.  So, for example $S_i$ will always be the $i$th strip of $\cB$, not the $i$th strip of $\cC$ which will be denoted $T_i$.  In the previous paragraph we saw that  $C$ starts with an $NI$ step.  Furthermore, the definition of $R'$ as a concatenation shows that $(k-t-1,0)$ is an $NI$ point of $C$.  So $\cC$ is again in case (a).  By induction $\io'(B')$ will be $B$ with its lowest row removed.  And the bottom row will be a left strip tiled by 
$$
(\io')^2(S_{t+1})\cdot R'\ce{t+1}^r = S_{t+1}\cdot R\ce{t+1}=R\fl{k-t-1}\cdot R\ce{t+1}=R.
$$
Thus $\io^2(B)=B$.  Also, using induction,
$$
\io^2(\cS) = \io'(\cT) =((\io')^2(S_1),\dots,(\io')^2(S_t)) = \cS.
$$
Hence $\io^2(\cB)=\cB$ as we wished to prove.

For the remaining three cases, much of the demonstration of being well defined is similar to what was done in case (a).  So we will just mention any important points of difference.  In case (b), the fact that $(k-t-1,0)$ is an $NL$ point for $B$ implies that there is a domino between squares $k-t-1$ and $k-t$ in the bottom row of $B$.  In particular, this means that $(k-t,0)$ is an $NI$ point for $B$ and so it is possible to take the first $k-t$ squares of $R$ when forming $R'$.  Note also that by definition of the right strip $R^R$, the domino just mentioned will cover squares $n-k+t$ and $n-k+t+1$ in the bottom row of $C$.  Thus a path starting at $(n-k+t,0)$ will be forced west and so this will an $NL$ point of $C$.  Furthermore, the first component of $\cT$ is $\io'(S_t)\cdot R\ce{t}^r$, where by induction $\io'(S_t)$ has length $[(n-1)-k+t]-t=n-k-1$.  So this gives an $NI$ point of $T_1$ with coordinates $(n-k-1,0)=((n-k+t)-t-1,0)$ and thus $\cC$ is in case (c).

To see that we have an involution in case (b), we have just noted that for $\cB$ in this case we have $\cC=\io(\cB)$ is in case (c). 
As usual, $\io'(B')$ returns the top rows of $B$ to what they were.   As for the bottom row we have, by definition of case (c) and the fact that $\cC$ is of type $(n,n-k+t,t)$, that it is a left strip tiled by
$$
(R')^r T_1^r\fl{n-(n-k+t)-t} = (R\fl{k-t})(R\ce{t}\cdot\io'(S_t)^r) \fl{k} = R. 
$$
Finally, we have
$$
\io^2(\cS)=(\io'(T_2),\dots,\io'(T_{t+1}))=((\io')^2(S_1),\dots,(\io')^2(S_{t-1}),\io'(T_1\fl{(n-k+t)-t-1}))
$$
where
$$
T_1\fl{(n-k+t)-t-1}=(\io'(S_t) \cdot R\ce{t}^r)\fl{n-k-1}=\io'(S_t).
$$
So by induction $\io^2(\cS)=\cS$ in this case as well.

The proof if $\cB$ is in case (c) is similar to the one for case (b) which is its inverse so this part of the demonstration will be omitted.  
Finally we turn to case (d).  To prove that this case is well defined, one again checks that the dominoes which force $(k,0)$ to be an $NL$ point of $B$ and $(k-t-1,0)$ to be an $NL$ point of $S_1$ appear in $T_1$ and $C$, respectively, so that $((n-k+t)-t-1,0)=(n-k-1,0)$ is an $NL$ point of $T_1$  and $(n-k+t,0)$ is an $NL$ point of $C$.  One then uses this fact to show that applying $\io$ twice is the identity.  But no new ideas appear so we will leave these details to the reader.
\eprf

\section{Catalan and Fuss-Catalan numbers}
\label{cfc}

The well-known {\em Catalan numbers} are given by
$$
C_n =\frac{1}{n+1}\binom{2n}{n}
$$
for $n\ge0$.  So the Lucas analogue is
$$
C_{\{n\}} = \frac{1}{\{n+1\}}{2n\brace n}.
$$
In 2010, Lou Shapiro suggested this definition.  Further, he asked whether this was a polynomial in $s$ and $t$ and, if so, whether it had a combinatorial interpretation.  There is a simple relation between $C_{\{n\}}$ and the Lucasnomials which shows that the answer to the first question is yes.  
This was first pointed out by Shalosh Ekhad~\cite{ekh:ssl}.
We will prove this equation combinatoriallly below.  We can now show that the second question also has an affirmative answer.
\bth
\label{Lcat:ptn}
Given $n\ge0$ there is a partition of $\cT(\de_{2n})$ such that $\{n\}!\{n+1\}!$ divides $\wt\beta$ for every block $\beta$.
\eth
\bprf
Given $T\in\cT(\de_{2n})$ we find the block  containing it as follows.  First construct a lattice path $p$ starting at $(n-1,0)$ and ending at $(2n,0)$ in exactly the same was as in the proof of Theorem~\ref{Lnom:ptn}.  Now put a tiling in the same block as $T$ if it agrees with $T$ on the left side of $NI$ steps and on the right side of $NL$ steps in all rows above the first row.  In the first row, the tiling on both sides of $p$ is arbitrary except for the required domino if $p$ begins with a $W$ step.  Since $p$ goes from $(n-1,0)$ to $(2n,0)$, the parts of the tiling which vary as in the Lucasnomial case contribute $\{n-1\}!\{n+1\}!$ to $\wt\beta$.  So we just need to show that the extra varying portion in the first row will give a factor of $\{n\}$.  If $p$ begins with an $N$ step, then the extra factor comes from the $n-1$ boxes to the left of this step which yields $\{n\}$.  If $p$ begins with $WN$, then this factor comes from the $n-1$ boxes to the right of the domino causing this $NL$ step, which again gives the desired $\{n\}$.
\eprf

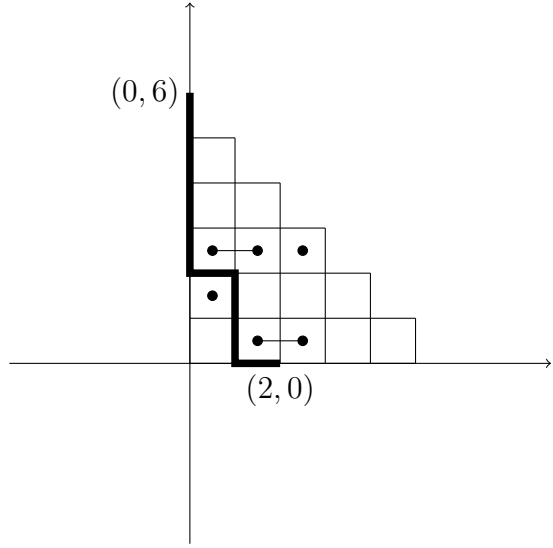
\begin{figure}
\bce
\begin{tikzpicture}[scale=0.6]
\draw (0,6)--(0,0);
\draw (1,5)--(1,0);
\draw (2,4)--(2,0);
\draw (3,3)--(3,0);
\draw (4,2)--(4,0);
\draw (5,1)--(5,0);
\draw (5,0)--(0,0);
\draw (5,1)--(0,1);
\draw (4,2)--(0,2);
\draw (3,3)--(0,3);
\draw (2,4)--(0,4);
\draw (1,5)--(0,5);
\draw[->] (0,-4)--(0,8);
\draw[->] (-4,0)--(8,0);
\draw[line width=1mm] (2,0)--(1,0)--(1,2)--(0,2)--(0,6);
\foreach \x in {1.5,2.5 } \filldraw (\x,.5) circle(3pt);
\foreach \x in {.5} \filldraw (\x,1.5) circle(3pt);
\foreach \x in {.5,1.5,2.5 } \filldraw (\x,2.5) circle(3pt);
\draw (1.5,.5)--(2.5,.5) (.5,2.5)--(1.5,2.5);
\node at (2,-0.6) {$(2,0)$};
\node at (-1,6) {$(0,6)$};
\end{tikzpicture}
\ece
\caption{A Catalan partial tiling}
\label{Cat:part}
\end{figure}

Again, we can associate with each block of the partition in the previous theorem a {\em Catalan partial tiling} which is like a binomial partial tiling except that the first row will be blank except for a domino if $p$ begins with a $W$ step.  We will sometimes omit the modifiers like ``binomial" and ``Catalan" if it is clear from context which type of partial tiling is intended.  Figure~\ref{Cat:part} illustrates a Catalan partial tiling

\bco
\label{Cpartial}
Given $n\ge0$ we have
$$
C_{\{n\}}=\sum_C \wt C
$$
where the sum is over all Catalan partial tilings $C$ associated with lattice paths from $(n-1,0)$ to $(0,2n)$ in $\de_{2n}$.
Thus $C_{\{n\}}\in\bbN[s,t]$.
\hqed
\eco

We can now give a combinatorial proof of the identity relating the Lucas-Catalan polynomials $C_{\{n\}}$ and the Lucasnomials which we mentioned earlier.
\bpr
\label{Csum}
For $n\ge2$ we have
$$
C_{\{n\}}={2n-1\brace n-1} + t{2n-1 \brace n-2}.
$$
\epr
\bprf
By Corollary~\ref{Cpartial}, it suffices to partition the Catalan partial tilings $P$ into two subsets whose weight generating functions are the two terms in the sum.  First consider the partial tilings associated with lattice paths $p$ whose first step is $N$.  Then the bottom row of $P$ is blank.  And the portion of $p$ in the remaining rows goes from $(n-1,1)$ to $(0,2n)$ inside $\de_{2n-1}$.  Thus the contribution of these partial tilings is ${2n-1\brace n-1}$.  If instead $p$ begins with $WN$, then there is a single domino in the first row which contributes $t$.  The rest of the path goes from $(n-2,1)$ to $(0,2n)$ inside  $\de_{2n-1}$ and so contributes 
${2n-1 \brace n-2}$ as desired.
\eprf

Note that this proposition is a Lucas analogue of the well-known identity
$$
C_n = \binom{2n-1}{n-1}-\binom{2n-1}{n-2}
$$
obtained when  $s=2$ and $t=-1$.

We now wish to study the Lucas analogue of the {\em Fuss-Catalan numbers} which are
$$
C_{n,k} =\frac{1}{kn+1}\binom{(k+1)n}{n}
$$
for $n\ge0$ and $k\ge1$.  Clearly $C_{n,1}=C_n$.  Consider the Lucas analogue
$$
C_{\{n,k\}}=\frac{1}{\{kn+1\}}{(k+1)n\brace n}.
$$

To prove the next result, it will be convenient to give coordinates to the squares of a Young diagram $\la$.  We will use brackets for these coordinates to distinguish them from the Cartesian coordinates we have been using for lattice paths.  Let $[i,j]$ denote the square in row $i$ from the bottom and column $j$ from the left.  Alternatively, if a square has northeast corner with Cartesian coordinates $(j,i)$ then the square's coordinates are $[i,j]$.
\bth
Given $n\ge0, k\ge1$ there is a partition of $\cT(\de_{(k+1)n})$ such that $\{n\}!\{kn+1\}!$ divides $\wt\beta$ for every block $\beta$.
\eth
\bprf
To find the block containing a tiling $T$ of $\de_{(k+1)n}$ we proceed as follows.  Consider the usual lattice path $p$ in $T$ starting at $(n-1,0)$ and ending at $((k+1)n,0)$.  If $p$ starts with an $N$ step, then we construct $\beta$ exactly as in the proof of Theorem~\ref{Lcat:ptn}.  In this case, the parts of the tiling which vary as in the Lucasnomial case contribute $\{n-1\}!\{kn+1\}!$ and the squares in the first row to the left of the $NI$ step give a factor of $\{n\}$ so we are done for such paths.

Now suppose $p$ begins $WN$.  It follows that there is a domino of $T$  between squares $[1,n-1]$ and $[1,n]$.  Also, there is no domino between squares $[1,(k+1)n-1]$ and $[1,(k+1)n]$ because the latter square is not part of $\de_{(k+1)n}$.  So there is a smallest index $m$ such that there is a domino between $[1,mn-1]$ and $[1,mn]$ but no domino between $[1,(m+1)n-1]$ and 
$[1,(m+1)n]$.  The block of $\beta$ will consist of all tilings agreeing with $T$ as for Lucasnomials in rows above the first.  And in the first row they agree with $T$ to the right of the $NL$ step except in the squares from $[1,mn+1]$ through $[1,(m+1)n-1]$ where the tiling is allowed to vary.  As in the previous paragraph, the variable parts of $\beta$ which are the same as for Lucasnomials contribute 
 $\{n-1\}!\{kn+1\}!$ while the variable portion to the right of the first $NL$ gives a factor of $\{n\}$.  This finishes the demonstration.
\eprf

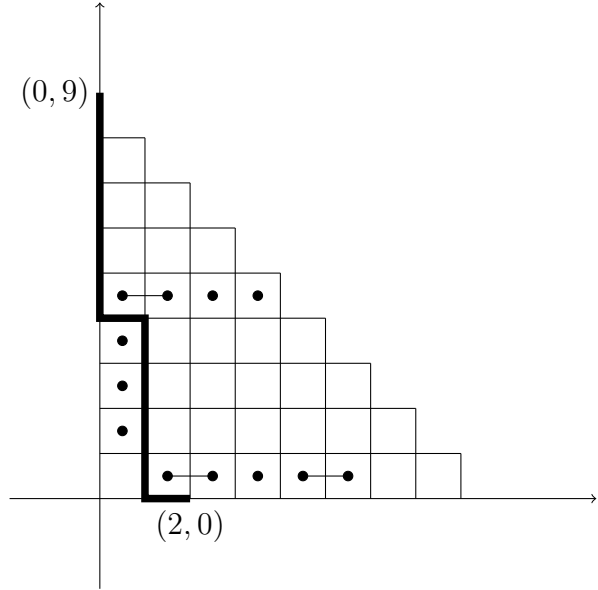
\begin{figure}
\bce
\begin{tikzpicture}[scale=0.6]
\draw (1,8)--(1,0) (2,7)--(2,0) (3,6)--(3,0) (4,5)--(4,0) (5,4)--(5,0) (6,3)--(6,0) (7,2)--(7,0) (8,1)--(8,0);
\draw (8,1)--(0,1) (7,2)--(0,2) (6,3)--(0,3) (5,4)--(0,4) (4,5)--(0,5) (3,6)--(0,6) (2,7)--(0,7) (1,8)--(0,8);
\draw[->] (0,-2)--(0,11);
\draw[->] (-2,0)--(11,0);
\draw[line width=1mm] (2,0)--(1,0)--(1,4)--(0,4)--(0,9);
\foreach \x in {1.5,2.5,...,5.5 } \filldraw (\x,.5) circle(3pt);
\foreach \x in {.5} \filldraw (\x,1.5) circle(3pt);
\foreach \x in {.5} \filldraw (\x,2.5) circle(3pt);
\foreach \x in {.5} \filldraw (\x,3.5) circle(3pt);
\foreach \x in {.5,1.5,...,3.5} \filldraw (\x,4.5) circle(3pt);
\draw (1.5,.5)--(2.5,.5) (4.5,.5)--(5.5,.5) (.5,4.5)--(1.5,4.5);
\node at (2,-0.6) {$(2,0)$};
\node at (-1,9) {$(0,9)$};
\end{tikzpicture}
\ece
\caption{A Fuss-Catalan partial tiling}
\label{FCat:part}
\end{figure}

As usual, we can represent a block $\beta$ of this partition by a {\em Fuss-Catalan partial tiling}.  Figure~\ref{FCat:part} displays such a tiling when $n=3$, $k=2$, and $m=2$ (in the notation of the previous proof).
\bco
\label{FCpartial}
Given $n\ge0, k\ge1$ we have
$$
C_{\{n,k\}}=\sum_P \wt P
$$
where the sum is over all Fuss-Catalan partial tilings associated with lattice paths going from $(n-1,0)$ to $(0,(k+1)n)$ in $\de_{(k+1)n}$.
Thus $C_{\{n,k\}}\in\bbN[s,t]$.
\hqed
\eco

The next result is proved in much the same way as Proposition~\ref{Csum} and so the demonstration is left to the reader.
\bpr
For $n\ge2, k\ge1$ we have

\vs{5pt}

\eqqed{
C_{\{n,k\}} = {(k+1)n-1\brace n-1} + \sum_{m=1}^k t^m \{n\}^{m-1} \{k-m+1\} {(k+1)n-1 \brace n-2}.
}
\epr

\section{Coxeter groups and $d$-divisible diagrams}
\label{cgd}

There is a way to associate a Catalan number and Fuss-Catalan numbers with any finite irreducible Coxeter group $W$.  This has lead to the area of reseach called Catalan combinatorics.  For more details, see the memoir of Armstrong~\cite{arm:gnp}.  The purpose of this section is to prove that for any $W$, the Coxeter-Catalan number is in $\bbN[s,t]$.  In fact, we will prove more general results using $d$-divisible Young diagrams.  This will also permit us to prove that for the infinite families of Coxeter groups, the Lucas-Fuss-Catalan analogue is in $\bbN[s,t]$.

\begin{figure}
$$
\barr{|c|l|l|}
\hline
W & d_1, \dots, d_n  & h \\
\hline
A_n & 2, 3, 4,\dots, n+1 & n+1\\
B_n & 2, 4, 6, \dots, 2n & 2n \\
D_n & 2, 4, 6, \dots, 2(n-1), n & 2(n-1) \hs{5pt} (\text{for $n\ge3$})\\
I_2(m) & 2,m & m \hs{5pt} (\text{for $m\ge2$})\\
H_3 & 2, 6,10 & 10\\
H_4 & 2, 12, 20, 30 & 30\\
F_4 & 2, 6, 8, 12 & 12\\
E_6 & 2, 5, 6, 8, 9, 12 & 12\\
E_7 & 2, 6, 8, 10, 12, 14, 18 & 18\\
E_8 & 2, 8, 12, 14, 18, 20, 24, 30 & 30\\
\hline
\earr
$$
\caption{finite irreducible Coxeter group degrees}
\label{deg}
\end{figure}

The finite Coxeter groups $W$ are those which can be generated by reflections.  Those which are irreducible have a well-known classification with four infinite families ($A_n$, $B_n$, $D_n$, and $I_2(m)$) as well as 6 exceptional groups ($H_3$, $H_4$, $F_4$, $E_6$, $E_7$, and $E_8$) where the subscript denotes the dimension $n$  of the space on which the group acts.  Associated with each finite irreducilble group is a set of {\em degrees} which are the degrees $d_1,\dots,d_n$ of certain polynomial invariants of the group.  The degrees of the various groups are listed in Figure~\ref{deg}.  The {\em Coxeter number} of $W$ is the largest degree and is denoted $h$.  One can now define the {\em Coxeter-Catalan number} of $W$ to be
$$
\Cat W = \prod_{i=1}^n \frac{h+d_i}{d_i}
$$
with corresponding {\em Lucas-Coxeter analogue}
$$
\Cat \{W\} = \prod_{i=1}^n \frac{\{h+d_i\}}{\{d_i\}}.
$$
If $W$ is of type $J_n$ for some $J$ then we will also use the notation  $J_{\{n\}}$ for $\{W\}$.
Directly from the definitions, $\Cat A_{n-1}=C_n$.  Also, after cancelling powers of $2$, we have 
$\Cat B_n=\binom{2n}{n}$.  But $\{2n\}\neq\{2\}\{n\}$ so we will have to find another way to deal with 
$\Cat B_{\{n\}}$.  In fact, we will be able to give a combinatorial interpretation when the numerator and denominator are both  constructed using ``Lucastorials" containing the integers divisible by some fixed integer $d\ge1$.

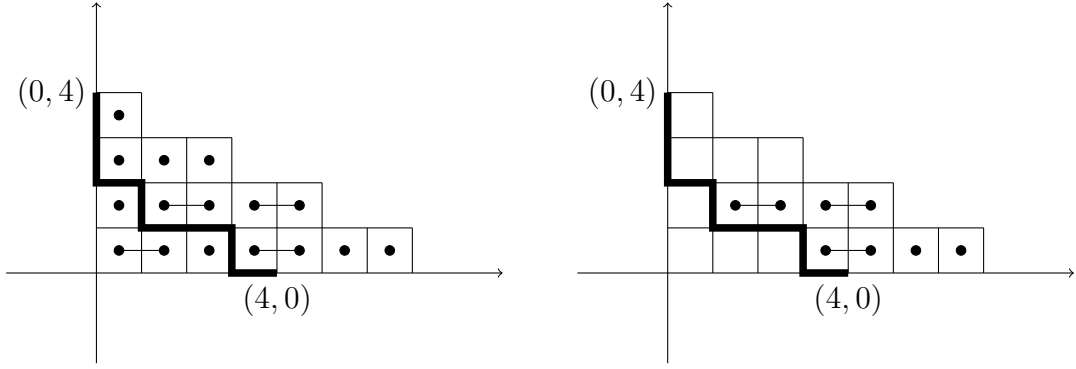
\begin{figure}
\bce
\begin{tikzpicture}[scale=0.6]
\draw (1,4)--(1,0) (2,3)--(2,0) (3,3)--(3,0) (4,2)--(4,0) (5,2)--(5,0) (6,1)--(6,0) (7,1)--(7,0);
\draw (7,1)--(0,1) (5,2)--(0,2) (3,3)--(0,3) (1,4)--(0,4);
\draw[->] (0,-2)--(0,6);
\draw[->] (-2,0)--(9,0);
\draw[line width=1mm] (4,0)--(3,0)--(3,1)--(1,1)--(1,2)--(0,2)--(0,4);
\foreach \x in {.5,1.5,...,6.5} \filldraw (\x,.5) circle(3pt);
\foreach \x in {.5,1.5,...,4.5} \filldraw (\x,1.5) circle(3pt);
\foreach \x in {.5,1.5,...,2.5} \filldraw (\x,2.5) circle(3pt);
\foreach \x in {.5} \filldraw (\x,3.5) circle(3pt);
\draw (0.5,.5)--(1.5,.5) (3.5,.5)--(4.5,.5) (1.5,1.5)--(2.5,1.5)  (3.5,1.5)--(4.5,1.5);
\node at (4,-0.6) {$(4,0)$};
\node at (-1,4) {$(0,4)$};
\end{tikzpicture}
\hs{20pt}
\begin{tikzpicture}[scale=0.6]
\draw (1,4)--(1,0) (2,3)--(2,0) (3,3)--(3,0) (4,2)--(4,0) (5,2)--(5,0) (6,1)--(6,0) (7,1)--(7,0);
\draw (7,1)--(0,1) (5,2)--(0,2) (3,3)--(0,3) (1,4)--(0,4);
\draw[->] (0,-2)--(0,6);
\draw[->] (-2,0)--(9,0);
\draw[line width=1mm] (4,0)--(3,0)--(3,1)--(1,1)--(1,2)--(0,2)--(0,4);
\foreach \x in {3.5,4.5,5.5,6.5} \filldraw (\x,.5) circle(3pt);
\foreach \x in {1.5,2.5,...,4.5} \filldraw (\x,1.5) circle(3pt);
\draw  (3.5,.5)--(4.5,.5)  (1.5,1.5)--(2.5,1.5)  (3.5,1.5)--(4.5,1.5);
\node at (4,-0.6) {$(4,0)$};
\node at (-1,4) {$(0,4)$};
\end{tikzpicture}
\ece
\caption{A $2$-divisible path on the left and corresponding  partial tiling on the right}
\label{ddiv:part}
\end{figure}

Define the {\em $d$-divisible Lucastorial} as
$$
\{n:d\}! = \{d\} \{2d\}\dots \{nd\}
$$
with corresponding {\em $d$-divisible Lucasnomial}
$$
{n:d\brace k:d} = \frac{\{n:d\}!}{\{k:d\}! \{n-k:d\}!}
$$
for $0\le k\le n$.
So we have
$$
\Cat  B_{\{n\}} = {2n:2\brace n:2}.
$$
Also define the {\em $d$-divisible staircase parttion} 
$$
\de_{n:d}  = (nd-1,(n-1)d-1,\dots,2d-1,d-1).
$$
The fact that $\wt\de_{n:d}=\{n:d\}!$ follows immediately from the definitions.
\bth
\label{dLnom:ptn}
Given $d\ge1$ and $0\le k\le n$ there is a partition of $\cT(\de_{n:d})$ such that $\{k:d\}!\{n-k:d\}!$ divides $\wt\beta$ for every block $\beta$.
\eth
\bprf
We determine the block $\beta$ containing a tiling $T$ by constructing a path $p$ from $(kd,0)$ to $(0,n)$ as follows.  The path takes an $N$ step if and only if three conditions are satisfied:  the two for Lucasnomial paths (the step does not cross a domino and stays within the Young diagram) together with the requirement that the $x$-coordinate of the $N$ step must be congruent to $0$ or $-1$ modulo $d$ with at most one $N$ step on each line of the latter type.  So $p$ starts by either goes north along $x=kd$ or, if there is a blocking domino, taking a $W$ step and going north along $x=kd-1$.  
In the first case it can take another $N$ step if not blocked, or go $W$ and then $N$ if it is.  In the second case, 
$p$  proceeds using $W$ steps to $((k-1)d,1)$ and either goes north from that lattice point or, if blocked, takes one more $W$ step to go north from $((k-1)d-1,1)$, etc.  See Figure~\ref{ddiv:part} for an example.  Call an $N$ step an $NI$ step if it has $x$-coordinate divisible by $d$ and an $NL$ step otherwise.  We now construct $\beta$ as for Lucasnomials: agreeing with $T$ to the left of $NI$ steps and to the right of $NL$ steps.  It is an easy matter to check that $\{k:d\}!\{n-k:d\}!$ is a factor of $\wt\beta$.
\eprf

The definition of {\em $d$-divisible partial tiling} (illustrated in Figure~\ref{ddiv:part}) and the next result are as expected.
\bco
\label{ddpartial}
Given  $d\ge1$ and $0\le k\le n$ we have
$$
{n:d\brace k:d}=\sum_P \wt P
$$
where the sum is over all $d$-divisible partial tilings associated with lattice paths going from $(kd,0)$ to $(0,n)$ in $\de_{n:d}$.
Thus ${n:d\brace k:d}\in\bbN[s,t]$.
\hqed
\eco  

The Lucas-Coxeter analogue for $D_n$ is
\beq
\label{CatD}
\Cat D_{\{n\}} = \frac{\{3n-2\}}{\{n\}}{2(n-1) : 2 \brace n-1 : 2}.
\eeq
Again, we will be able to prove a $d$-divisible generalization of this result.  But first we need a result of Hoggatt and Long~\cite{hl:dpg} about the divisibility of polynonials in the Lucas sequence.  (In their paper they only prove the divisbility statement, but the fact that the quotient is in $\bbN[s,t]$ follows easily from their demonstration.)

\begin{figure}[t]
\begin{tikzpicture}[scale=0.6]
\draw (1,9)--(1,1) (2,8)--(2,1) (3,8)--(3,1) (4,7)--(4,1) (5,7)--(5,0) (6,6)--(6,0) (7,6)--(7,0) 
   (8,5)--(8,0) (9,5)--(9,0) (10,4)--(10,0) (11,4)--(11,0) (12,3)--(12,0) (13,3)--(13,0) (14,2)--(14,0) (15,2)--(15,0) (16,1)--(16,0) (17,1)--(17,0);
\draw (17,1)--(0,1) (15,2)--(0,2) (13,3)--(0,3) (11,4)--(0,4) (9,5)--(0,5) (7,6)--(0,6) (5,7)--(0,7) (3,8)--(0,8) (1,9)--(0,9);
\draw[->] (0,-1)--(0,10);
\draw[->] (-1,0)--(18,0);
\draw[line width=1mm] (10,0)--(10,1)--(10,2)--(9,2)--(9,3)--(8,3)--(7,3)--(7,4)--(6,4)--(6,5)--(6,6)--(5,6)--(5,7)--(3,7)--(3,8)--(1,8)--(1,9)--(0,9);
\foreach \x in {5.5,6.5,...,16.5} \filldraw (\x,.5) circle(3pt);
\foreach \x in {.5,1.5,...,14.5} \filldraw (\x,1.5) circle(3pt);
\foreach \x in {.5,1.5,...,12.5} \filldraw (\x,2.5) circle(3pt);
\foreach \x in {.5,1.5,...,10.5} \filldraw (\x,3.5) circle(3pt);
\foreach \x in {.5,1.5,...,8.5} \filldraw (\x,4.5) circle(3pt);
\foreach \x in {.5,1.5,...,6.5} \filldraw (\x,5.5) circle(3pt);
\foreach \x in {.5,1.5,...,4.5} \filldraw (\x,6.5) circle(3pt);
\foreach \x in {.5,1.5,...,2.5} \filldraw (\x,7.5) circle(3pt);
\foreach \x in {.5} \filldraw (\x,8.5) circle(3pt);
\draw (7.5,2.5)--(8.5,2.5) (9.5,2.5)--(10.5,2.5) (2.5,2.5)--(3.5,2.5)  (3.5,3.5)--(4.5,3.5) (7.5,3.5)--(8.5,3.5)
         (9.5,.5)--(8.5,.5) (14.5,.5)--(15.5,.5) (8.5,1.5)--(9.5,1.5)  (3.5,1.5)--(4.5,1.5)
         (.5,1.5)--(1.5,1.5) (11.5,1.5)--(12.5,1.5) (.5,4.5)--(1.5,4.5) (4.5,4.5)--(5.5,4.5)  (1.5, 5.5)--(2.5,5.5) (6.5,4.5)--(7.5,4.5)
         (1.5,5.5)--(2.5,5.5) (.5,7.5)--(1.5,7.5)  (1.5,6.5)--(2.5,6.5);
\node at (10,-0.6) {$(10,0)$};
\node at (-1,9) {$(0,9)$};
\end{tikzpicture}

\hs{10pt}

\begin{tikzpicture}[scale=0.6]
\draw (1,9)--(1,1) (2,8)--(2,1) (3,8)--(3,1) (4,7)--(4,1) (5,7)--(5,0) (6,6)--(6,0) (7,6)--(7,0) 
   (8,5)--(8,0) (9,5)--(9,0) (10,4)--(10,0) (11,4)--(11,0) (12,3)--(12,0) (13,3)--(13,0) (14,2)--(14,0) (15,2)--(15,0) (16,1)--(16,0) (17,1)--(17,0);
\draw (17,1)--(0,1) (15,2)--(0,2) (13,3)--(0,3) (11,4)--(0,4) (9,5)--(0,5) (7,6)--(0,6) (5,7)--(0,7) (3,8)--(0,8) (1,9)--(0,9);
\draw[->] (0,-1)--(0,10);
\draw[->] (-1,0)--(18,0);
\draw[line width=1mm] (10,0)--(10,1)--(10,2)--(9,2)--(9,3)--(8,3)--(7,3)--(7,4)--(6,4)--(6,5)--(6,6)--(5,6)--(5,7)--(3,7)--(3,8)--(1,8)--(1,9)--(0,9);
\foreach \x in {5.5,6.5,...,9.5} \filldraw (\x,.5) circle(3pt);
\foreach \x in {.5,1.5,...,9.5} \filldraw (\x,1.5) circle(3pt);
\foreach \x in {9.5,10.5,...,12.5} \filldraw (\x,2.5) circle(3pt);
\foreach \x in {7.5,8.5,...,10.5} \filldraw (\x,3.5) circle(3pt);
\foreach \x in {.5,1.5,...,5.5} \filldraw (\x,4.5) circle(3pt);
\foreach \x in {.5,1.5,...,5.5} \filldraw (\x,5.5) circle(3pt);
\draw (9.5,2.5)--(10.5,2.5)   (7.5,3.5)--(8.5,3.5)
         (9.5,.5)--(8.5,.5) (8.5,1.5)--(9.5,1.5)  (3.5,1.5)--(4.5,1.5)
         (.5,1.5)--(1.5,1.5) (.5,4.5)--(1.5,4.5) (4.5,4.5)--(5.5,4.5)  (1.5, 5.5)--(2.5,5.5) 
         (1.5,5.5)--(2.5,5.5) ;
\node at (10,-0.6) {$(10,0)$};
\node at (-1,9) {$(0,9)$};
\end{tikzpicture}
\capt{A lattice path and partial tiling for $D_{\{5\}}$ \label{D5}}
\end{figure}

\bth[\cite{hl:dpg}]
\label{hl}
For positive integers $m,n$ we have $m$ divides $n$ if and only if $\{m\}$ divides $\{n\}$.
In this case $\{m\}/\{n\}\in\bbN[s,t]$.
\hqed.
\eth

If $\la,\mu$ are Young diagrams with $\mu\sbe\la$, then the corresponding {\em skew diagram}, $\la/\mu$, consists of all the boxes in $\la$ but not in $\mu$.  The skew diagram used in Figure~\ref{D5} is $\de_{9:2}/(5)$.  We will use skew diagrams to prove the following result which yields~\ree{CatD} as a special case.
\bth
Given $d\ge1$ we have
\beq
\label{dCatD}
\frac{\{(d+1)n-d\}}{\{n\}}{2(n-1) : d \brace n-1: d}
\eeq
is in $\bbN[s,t]$.
\eth
\bprf
Tile the rows of the skew shape $\de_{2n-1:d}/(d-1)n$. It is easy to see that the corresponding generating function is the numerator of equation~\ree{dCatD} where the bottom row contributes $\{(2n-1)d-(d-1)n\}=\{(d+1)n-d\}$ which is  the numerator of the fractional factor.  As usual, we group the tilings into blocks $\be$ and show that the denominator of~\ree{dCatD} divides the weight of each block.  Given a tiling $T$ we find its block by starting a lattice path $p$ at $(nd,0)$ and using exactly the same rules as in the proof of Theorem~\ref{dLnom:ptn}.  See the upper diagram in Figure~\ref{D5} for an example when $d=2$ and $n=5$.  We now let the strips to the side of each north step either be fixed or vary, again as dictated in Theorem~\ref{dLnom:ptn}'s demonstration.  The bottom diagram in Figure~\ref{D5} shows the partial tiling corresponding to the upper diagram.  There are two cases.

If $p$ starts with an $N$ step then its right side contributes $\{(2n-1)d-nd\}=\{(n-1)d\}$.  Now $p$ enters the top $2n-2$ rows which form a $\de_{2n-2:d}$ at $x$-coordinate $nd$.  It follows that the contribution of this portion of $p$ to $\wt\be$ is 
$\{n:d\}!\{n-2:d\}!$.  So the total contribution to the weight of the variable parts of each row is
$$
\{(n-1)d\}\cdot \{n:d\}!\{n-2:d\}! = \{nd\}\cdot  \{n-1:d\}!\{n-1:d\}!.
$$
Thanks to Theorem~\ref{hl}, this is divisible by $\{n\}\cdot  \{n-1:d\}!\{n-1:d\}!$ which is the denominator of~\ree{dCatD} so we are done with this case.

If $p$ starts with $WN$, then the left side of $p$ gives a contribution of $\{nd-n(d-1)\}=\{n\}$.  Because of the rule that $p$ can take at most one step on a vertical line of the form $=kd-1$, the path must now continue to take west steps until it reaches $((n-1)d,1)$.  Now it can enter the upper rows of the diagram from this point contributing a factor of $ \{n-1:d\}!\{n-1:d\}!$ to $\wt\be$.  Multiplying the two contributions gives exactly the denominator of~\ree{dCatD} which finishes this case and the proof.
\eprf

We note that replacing $n$ by $n+1$ in equation~\ree{dCatD} we obtain
$$
\frac{\{(d+1)n+1\}}{\{n+1\}}{2n : d \brace n: d}
$$
which, for $d=1$, is just a multiple of the $n$th Lucas-Catalan polynomial.  We also note that one can generalize even further. 
Given positive integers satisfying  $l<kd<md$, we consider starting lattice paths from $(kdn,0)$ in the skew diagram 
$\de_{m(n-1)+1:d}/(ld)$ using the same rules as in the previous two proofs.  This yields the following result whose demonstration is similar enough to those just given that we omit it.  But note that we get the previous theorem as the special case $l=d-1$, $k=1$, and $m=2$.
\bth
\label{genCatD}
Given positive integers satisfying $l<kd<md$ we have 
$$
\frac{\{(dm-l)n-(m-1)d\}}{\{gn\}} {m(n-1) : d\brace kn-1 :d}
$$
is in $\bbN[s,t]$ where $g=\gcd(kd,kd-l)$.\hqed
\eth

We finally come to our main theorem for this section.
\bth
\label{CoxCat}
If $W$ is a finite irreducible Coxeter group $W$ then $\Cat \{W\}\in\bbN[s,t]$.
\eth
\bprf
We have already proved the result in types $A$, $B$, and $D$.  And for the exceptional Coxeter groups, we have verified the claim by computer.  So we just need to show that
$$
\Cat I_{\{2\}}(m) =\frac{\{m+2\}\{2m\}}{\{2\}\{m\}}
$$
is in $\bbN[s,t]$.  But this follows easily from Theorem~\ref{hl}.  Indeed,  if $m$ is odd then the  relatively prime polynomials $\{2\}=s$ and $\{m\}$ both divide $\{2m\}$.  It follows that the same is true of their product which completes this case.  If $m$ is even then $\{2\}$ divides $\{m+2\}$ and $\{m\}$ divides $\{2m\}$.  So, again, we have a polynomial quotient and all quotients have nonnegative integer coefficients.
\eprf

We now turn to the Fuss-Catalan case.  For any finite Coxeter group $W$ and positive integer $k$ there is a {\em Coxeter-Fuss-Catalan number} defined by
$$
\Cat^{(k)} W = \prod_{i=1}^n \frac{kh + d_i}{d_i}.
$$
So, there is a corresponding {\em Lucas-Coxeter-Fuss-Catalan analogue} given by
$$
\Cat^{(k)} \{W\} = \prod_{i=1}^n \frac{\{kh + d_i\}}{\{d_i\}}.
$$
In particular, $\Cat^{(k)} A_{\{n-1\}} = C_{\{n,k\}}$.

\bth
If $W=A_n, B_n, D_n, I_2(m)$ then  $\Cat^{(k)} \{W\}\in\bbN[s,t]$.
\eth
\bprf
We have already shown this for $A_n$ in Corollary~\ref{FCpartial}.
For type $B$ we find that
$$
\Cat^{(k)} B_{\{n\}} = \prod_{i=1}^n \frac{\{2kn + 2i\}}{\{2i\}} = {(k+1)n :2 \brace n:2}
$$
which is a polynomial in $s$ and $t$ with nonnegative coefficients by Corollary~\ref{ddpartial}.
In the case of type $D$ we see 
\begin{align*}
\Cat^{(k)} D_{\{n\}} &= \frac{\{2k(n-1) + n\}}{\{n\}} \prod_{i=1}^{n-1} \frac{\{2k(n-1) + 2i\}}{\{2i\}}\\[10pt]
&= \frac{\{(2k+1)n - 2k\}}{\{n\}} {(k+1)(n-1) : 2 \brace n-1 : 2}
\end{align*}
which is in $\bbN[s,t]$ by Theorem~\ref{genCatD}.

For $\Cat^{(k)} I_{\{2\}}(m)$ we can use an argument similar to that of the proof of Theorem~\ref{CoxCat}.
We have that
$$
\Cat^{(k)} I_{\{2\}}(m) = \frac{ \{km + 2\} \{(k+1)m \}}{\{2\} \{m\}}
$$
and can consider the parity of $m$ and $k$.
If $m$ or $k$ is even then  $\{2\}$ divides $\{km+2\}$ and $\{m\}$ divides $\{(k+1)m\}$.
If $m$ and $k$ are both odd then  $\{2\}$ and $\{m\}$ are relatively prime and both divide $\{(k+1)m\}$.
And in all cases the quotients have coefficients in $\bbN$.
\eprf

\section{Comments and future work}

Here we will collect various observations and open problems in the hopes that the reader will be tempted to continue our work.

\subsection{Coefficients}  

Note that we can write
$$
\{n\} = \sum_k a_k s^{n-2k-1} t^k
$$
where the $a_k$ are positive integers and $0\le k\le (n-1)/2$.  We call $a_0,a_1,\dots$ the {\em coefficient sequence} of $\{n\}$ and note that any of our Lucas analogues considered previously will also correspond to such a sequence.  There are several properties of sequences of real numbers which are common in combinatorics, algebra, and geometry.  One is that the sequence is {\em unimodal} which means that there is an index $m$ such that 
$$
a_0\le a_1\le \dots \le a_m\ge a_{m+1}\ge \dots.
$$
Another is that the sequence is {\em log concave} which is defined by the inequality
$$
a_k^2 \ge a_{k-1} a_{k+1}
$$
for all $k$, where we assume $a_k=0$ if the subscript is outside of the range of the sequence.  Finally, we can consider the generating function 
$$
f(y) = \sum_{k\ge0} a_k y^k
$$
and ask for properties of its roots.  For more information about such matters, see the survey articles of Stanley~\cite{sta:lus} and Brenti~\cite{bre:lus}.  In particular, the following result is well known and straightforward to prove.
\bpr
If $a_0,a_1,\dots$ is a sequence of positive reals then its generating function having real roots implies that it is log concave.  And if the sequence is log concave then it is unimodal.\hqed
\epr

To see which of these properties are enjoyed by our Lucas analogues, it will be convenient to make a connection with Chebyshev polynomials.  The {\em Chebyshev polynomials of the second kind}, $U_n(x)$, are defined recursively by $U_0(x)=1$, $U_1(x)=2x$, and for $n\ge 2$
$$
U_n(x) = 2x U_{n-1}(x) - U_{n-2}(x).
$$
It follows immediately that 
\beq
\label{Un}
\{n\} = U_{n-1}(x) 
\eeq
if we set $s=2x$ and $t=-1$.
\bth
If the Lucas analogue of a quotient of products is a polynomial then it has a coefficient generating function which is real rooted.  So if the coefficient sequence consists of  positive integers then the sequence is log concave and unimodal.
\eth
\bprf
From the previous proposition, it suffices to prove the first statement.
It is well known and easy to prove by using angle addition formulas that
$$
U_n(\cos\th) = \frac{\sin(n+1)\th}{\sin\th}.
$$
It follows that the roots of $U_n(x)$ are 
$$
x=\cos\frac{k\pi}{n+1}
$$
for $0<k<n+1$ and so real.

By equation~\ref{Un} we see that $\{n\}$ and $U_{n-1}(x/2)$ have the same coefficient sequence except that in the former all coefficients are positive and in the latter signs alternate.  Now if we take a quotient of products of the $\{n\}$ which is a polynomial  $p(s,t)$, then the corresponding quotient of products where $\{n\}$ is replaced by $U_{n-1}(x/2)$ will be a polynomial $q(x)$.  Further, from the paragraph above, $q(x)$ will have real roots.  It follows that  the coefficient sequence of $p(s,t)$ (which is obtained from $q(x)$ by removing zeros and making all coefficients positive) has a generating function with only real roots.
\eprf

\subsection{Rational Catalan numbers}

Rational Catalan numbers generalize the ordinary Catalans and have many interesting connections with Dyck paths, noncrossing partitions, associahedra, cores of integer partitions, parking functions, and Lie theory.  Let $a,b$ be positive  integers whose greatest common divisor is  $(a,b)=1$.  Then the associated {\em ratioinal Catalan number} is
$$
\Cat(a,b) = \frac{1}{a+b}\binom{a+b}{a}.
$$
In particular, it is easy to see that  $\Cat(n,n+1)=C_n$.  We note that the the relative prime condition is needed to ensure that $\Cat(a,b)$ is an integer.  As usual, consider
$$
\Cat\{a,b\}=\frac{1}{\{a+b\}}{a+b \brace a}.
$$
We will now present a proof that $\Cat\{a,b\}$ is a polynomial in $s,t$ which was obtained by the Fields Institute Algebraic Combinatorics Group~\cite{abcdl} for the $q$-Fibonacci analogue and works equally well in our context.  First we need the following lemma.
\ble[\cite{hl:dpg}]
We have $(\{m\},\{n\})=\{(m,n)\}$.\hqed
\ele

\bth[\cite{abcdl}]
If $(a,b)=1$ then $C\{a,b\}$ is a polynomial in $s,t$.
\eth
\bprf
Consider the quantity
$$
p = \frac{\{a+b\}!}{\{a-1\}!\{b\}!} = \{a+b\} {a+b-1 \brace a-1}.
$$
From the second expression for $p$, it is clearly a polynomial in $s$ and $t$.
Note that $\{a\}$ divides evenly into $p$ because $p/\{a\}= {a+b \brace a}$.
Similarly, $\{a+b\}$ divides $p$ since $p/\{a+b\}={a+b-1 \brace a-1}$.  But 
$(a,a+b)=1$ and so, by the lemma, $\{a\}\{a+b\}$ divides into $p$.
It follows that $p/(\{a\}\{a+b\}) = C\{a,b\}$ is a polynomial in $s,t$.
\eprf

Despite this result, we have been unable to find a combinatorial interpretation for $C\{a,b\}$ or prove that its coefficients are nonnegative integers, although this has been checked by computer for $a,b\le 50$.  

\subsection{Narayana numbers.}

For $1\le k\le n$ we define the {\em Narayana number}
$$
N_{n,k} = \frac{1}{n} \binom{n}{k} \binom{n}{k-1}.
$$
It is natural to consider the Narayana numbers in this context because $C_n=\sum_k N_{n,k}$.  Further discussion of Narayana numbers can be found in the paper of Br\"and\'en~\cite{bra:qnn}.  In our usual manner, let
$$
N_{\{n,k\}} =\frac{1}{\{n\}} {n \brace k} {n\brace k-1}.
$$
\bcon
For all $1\le k\le n$ we have $N_{\{n,k\}}\in\bbN[s,t]$.
\econ
This conjecture has been checked by computer for $n\le 100$.  One could also consider Narayana numbers for other Coxeter groups.

\subsection{Coxeter groups again}

As is often true in the literature on Coxeter groups, our proof of Theorem~\ref{CoxCat} is case by case.  It would be even better if a case-free demonstration could be found.  One could also hope for a closer connection between the geometric or algebraic properties of Coxeter groups and our combinatorial constructions.  Alternatively, it would be quite interesting to give a proof of Theorem~\ref{CoxCat} by weighting one of the standard objects counted by $\Cat W$ such as $W$-noncrossing partitions.

\vs{10pt}

{\bf Acknowledgment.}  We had helpful  discussions with Nantel Bergeron, Cesar Ceballos, and Volker Strehl about this work.



\nocite{*}
\bibliographystyle{alpha}

\newcommand{\etalchar}[1]{$^{#1}$}

\end{document}